\documentclass[11pt]{amsart}

\usepackage[utf8]{inputenc}
\usepackage{bm}
\usepackage{fontenc}
\usepackage{amsfonts}
\usepackage{amssymb}
\usepackage{amsmath}
\usepackage{amsthm}
\usepackage{mathtools}
\usepackage{enumerate}
\usepackage{enumitem}
\usepackage[pagebackref,colorlinks,linkcolor=blue,citecolor=blue,urlcolor=blue,hypertexnames=true]{hyperref}
\usepackage[pagebackref,hypertexnames=true]{hyperref}
\usepackage{enumitem}
\usepackage{mathrsfs}
\usepackage{tikz}
\usetikzlibrary{calc}
\usepackage{marginnote}
\usepackage{xcolor,enumitem}
\usepackage{soul}
\usepackage{tikz}
\usepackage[all,cmtip]{xy} 
\usepackage{caption}

\newcommand{\cstarx }[1]{\mathcal {#1}}
\newcommand{\N }{\mathbb {N}}
\newcommand{\C }{\mathbb {C}}
\newcommand{\eps }{\varepsilon }
\usepackage{bbm}
\newcommand{\SOTh }{\mathrm {SOT}\text {-}}
\newcommand{\cstu }{\mathrm {C}^*_u}
\newcommand{\cstql }{\mathrm {C}^*_{\textit {ql}}}
\newcommand{\cstar }{$\mathrm {C}^*$}
\newcommand{\cB }{\mathcal {B}}
\newcommand{\cK }{\mathcal {K}}

\newcommand{\B }{\cstarx {B}}

\newcommand{\R }{\mathbb {R}}
\newcommand{\cL }{\mathcal L}
\newcommand{\Adm}{\mathcal L _{adm}(\ell_2(X))}

\numberwithin{equation}{subsection}

\theoremstyle{theorem}
   \newtheorem{theorem}[equation]{Theorem}
   \newtheorem{lemma}[equation]{Lemma}
   \newtheorem{proposition}[equation]{Proposition}
   \newtheorem{corollary}[equation]{Corollary}

\theoremstyle{definition}
   \newtheorem{definition}[equation]{Definition}
   \newtheorem{claim}[equation]{Claim}
   \newtheorem{convention}[equation]{Convention}
   \newtheorem{assumption}[equation]{Assumption}
   \newtheorem{question}[equation]{Question}
   
\newtheorem*{AIusage}{A remark on this version of the manuscript and on AI usage}

\theoremstyle{remark}
   \newtheorem{remark}[equation]{Remark}
   \newtheorem{example}[equation]{Example}
\newtheorem*{acknowledgement*}{Acknowledgement}
 
\DeclareMathOperator{\dom }{dom}
\DeclareMathOperator{\ran }{ran}

\DeclareMathOperator{\propg }{prop}

\begin{document}

\title[Flows on uniform Roe algebras]{ Flows  on uniform Roe algebras}%

\date{\today } 
 
\author[B. M. Braga]{Bruno M. Braga}
\address[B. M. Braga]{IMPA, Estrada Dona Castorina 110, 22460-320, Rio de Janeiro, Brazil}
\email{demendoncabraga@gmail.com}
\urladdr{https://sites.google.com/view/demendoncabraga}
\thanks{B. M. Braga  was partially supported by FAPERJ, grant E-26/204.317/2025,  by CNPq, grant 303182/2026-1, and by Serrapilheira, grant R-2501-51476.}

\author[A. Buss]{Alcides Buss}
\address[A. Buss]{Universidade Federal de Santa Catarina, 88040-970 Florianopolis SC, Brazil}
\email{alcides.buss@ufsc.com}
\urladdr {http://www.mtm.ufsc.br/~alcides/}

\author[R. Exel]{Ruy Exel}
\address[R. Exel]{Universidade Federal de Santa Catarina, 88040-970 Florianopolis SC, Brazil}
\email{ruyexel@gmail.com}
\urladdr{http://www.mtm.ufsc.br/~exel/}
\thanks{A. Buss and R. Exel were partially supported by  CNPq, Brazil}
 
\maketitle

\begin{abstract}
For a uniformly locally finite metric space $(X, d)$, we investigate \emph{coarse} flows on its uniform Roe algebra $\cstu(X)$, defined as one-parameter groups of automorphisms whose differentiable elements include all partial isometries arising from partial translations on $X$. We first show that any flow $\sigma$ on $\cstu(X)$ corresponds to a (possibly unbounded) self-adjoint operator $h$ on $\ell_2(X)$ such that $\sigma_t(a) = e^{ith} a e^{-ith}$ for all $t \in \mathbb{R}$, allowing us to focus on operators $h$ that generate flows on $\cstu(X)$.

Assuming Yu's property A, we prove that a self-adjoint operator $h$ on $\ell_2(X)$ induces a coarse flow on $\cstu(X)$ if and only if $h$ can be expressed as $h = a + d$, where $a \in \cstu(X)$ and $d$ is a diagonal operator with entries forming a coarse function on $X$. We further study cocycle equivalence and cocycle perturbations of coarse flows, showing that, under property A, any coarse flow is a cocycle perturbation of a diagonal flow. Finally, for self-adjoint operators $h$ and $k$ that induce coarse flows on $\cstu(X)$, we characterize conditions under which the associated flows are either cocycle perturbations of each other or cocycle conjugate to each other. In particular, if $h - k$ is bounded, then the flow induced by $h$ is a cocycle perturbation of the flow induced by $k$.
\end{abstract}

\section{Introduction}
\setcounter{subsection}{1}

Given a (uniformly locally finite) metric space $X$, its uniform Roe algebra $\cstu (X)$ is a $\mathrm C^*$-algebra which captures the large scale geometry of $X$ in $\mathrm C^*$-algebraic terms (see Subsection \ref{SubsectionURA} for its precise definition). These algebras were introduced by J. Roe in the study of the index theory of elliptic operators on noncompact manifolds (\cite {Roe1988,Roe1993}) and are of interest to researchers working on the coarse Baum-Connes and Novikov conjectures (\cite {Yu2000}). Recently, mathematical physicists have started using uniform Roe algebras in the study of topological insulators (see \cite {Kubota2017,EwertMeyer2019,Jones2021CommMathPhys,LudewigThiang2021CommMathPhys,Bourne2022JPhys}).

Motivated by mathematical physics, \cite {BragaExel2023KMS} initiated the treatment of flows and KMS states on uniform Roe algebras. In a nutshell, the aforementioned article identified (and analysed) a naturally occurring class of flows on uniform Roe algebras: given a uniformly locally finite metric space $(X,d)$ and a map $h\colon X\to \R $, we can interpret $h$ as a multiplication (unbounded) operator on $\ell _2(X)$ --- the Hilbert space of square-summable functions $X\to \C $. Throughout these notes, we use the   convention that an \emph{operator} on a Hilbert space $H$ means a \emph{not necessarily bounded} operator defined on a vector \emph{subspace} of $H$.  A multiplication operator induced by such $h\colon X\to \R$  is clearly self-adjoint and diagonal with respect to the canonical orthonormal basis $(\delta _x)_{x\in X}$ of $\ell _2(X)$. In particular, $e^{ith}$ is a well-defined unitary operator on $\ell _2(X)$ which is also diagonal and the formula
  \begin{equation}
  \label{Eq11mar24FlowDiag}
  \sigma _{h,t}(a)=e^{ith}ae^{-ith},
  \end{equation}
  for $t\in \R $ and $a\in \cstu (X)$, defines a \emph{pre-flow}  (i.e., a not necessarily continuous action) $\sigma _h$ of $\R $ on $\cstu (X)$ by automorphisms. Moreover, the large-scale geometric properties of $h$  may be used to characterize when  $\sigma _h$ is strongly continuous, in which case we call it a \emph{flow}.  Recall that a map
  $h=(h_x)_{x\in X}\colon X\to \R $ is \emph{coarse} if, for all $r>0$, there is $s>0$, such that
  \[
  d(x,z)\leq r\ \text { implies }\ |h_x-h_z|\leq s
  \]
  for all $x,z\in X$.\footnote{The terms \emph{bornologous} and \emph{controlled} are often used in the literature instead of \emph{coarse}.} We then have the following:

  \begin{proposition}
  \emph{(}\cite [Proposition 2.1]{BragaExel2023KMS}\emph{)}. Let $X$ be a uniformly locally finite metric space and $\sigma _h$ be given by a map $h\colon X\to \R $ as in \eqref{Eq11mar24FlowDiag}. Then $\sigma _h$ is a flow on $\cstu (X)$ if and only if $h\colon X\to \R $ is coarse. \label{Prop2.1BragaExel}
  \end{proposition}

Flows of this sort are called \emph{diagonal flows}. The goal of the present paper is to extend the study of flows on uniform Roe algebras outside the class of diagonal flows. For that, we develop a theory of \emph{noncommutative coarse maps} or, equivalently, of \emph{coarse  operators}, and study  cocycle perturbation and cocycle conjugacy for flows on uniform Roe algebras.

One of the main reasons justifying one's interest in cocycle equivalence for flows is that, when two flows are cocycle conjugate to each other, then the fields of KMS states are isomorphic (\cite[Proposition 2.1]{Kishimoto2000RepMathPhy}).  So, in particular, all of the information about phase transitions may be passed from one flow to the other.  Since a lot is now known about KMS states for diagonal flows (\cite {BragaExel2023KMS}), the results herein extend that knowledge to a vastly larger class of flows.

\begin{acknowledgement*}
 We would like to express our thanks to Stuart White for some very inspiring and visionary conversations during a visit of the third named author to Oxford, when some conjectures we attempted to solve here were first sketched.
\end{acknowledgement*}

We now describe the main findings of this paper in more detail.

  \subsection{Uniform Roe algebras}\label{SubsectionURA} All metric spaces considered herein will be not only discrete but also uniformly locally finite.  Recall that a metric space $X$ is \emph{uniformly locally finite} (abbreviated as \emph{u.l.f}.\ from now on) if
  \[
  \sup _{x\in X}|B_r(x)|<\infty \ \text { for all }\ r>0,
  \] where $|B_r(x)|$ denotes the number of elements of the closed ball in $X$ of radius $r$ centered at $x$. In this case, $X$ is necessarily countable and the metric topology is discrete.

Given a u.l.f.\ metric space $(X,d)$, we identify $\ell _\infty (X)$ --- the $\mathrm C^*$-algebra of bounded, complex valued maps on $X$ --- with the algebra of diagonal operators on $\ell _2(X)$ in the canonical way. The uniform Roe algebra of $X$ is then the $\mathrm C^*$-subalgebra of $\cB (\ell _2(X))$\footnote {Given a Hilbert space $H$, $\cB (H)$ denotes the space of all bounded operators on $H$.}  generated by the partial isometries on $\ell _2(X)$ given by partial translations of $X$. Precisely, the notion of a translation function on a vector space is adapted to arbitrary metric spaces as follows: a partial bijection $f\colon \dom (f)\subseteq X\to \mathrm {ran}(f)\subseteq X$ is said to be a
  \emph{partial translation of $X$} if
  \[
  \sup _{x\in \dom (f)} d(x,f(x))<\infty .
  \]
  Any such $f$ induces a sort of \emph{partial shift} operator  $v_f$ on $\ell _2(X)$ determined by
  \[
  v_f\delta _x=\left \{
  \begin{array}
  {ll}
    \delta _{f(x)}, & \text { if } x\in \dom (f),\\
    0, & \text { if } x\not \in \dom (f).
  \end{array}
  \right .
  \]
  We can now make the definition of the uniform Roe algebra of a u.l.f.\ metric space precise:

  \begin{definition}\label{DefineCstu} The \emph{uniform Roe algebra} of a u.l.f.\ metric space $X$, denoted $\cstu (X)$, is the $\mathrm {C}^*$-subalgebra of $\cB (\ell _2(X))$ generated by all $v_f$, where $f$ is a partial translation of $X$.
  \end{definition}

  While we will only deal with u.l.f.\ spaces here, one can define uniform Roe algebras for more general spaces. This is, for instance, done in \cite{MR4411373,2024arXiv240313624M}.
  
If $A$ is any subset of $X$, observe that the orthogonal projection $p_A$ from $\ell _2(X)$ onto  $\ell _2(A)\subseteq \ell_2(X)$ coincides with the operator $v_f$, where $f$ is the identity function on $\dom (f) \coloneqq A$.  As the set of all such operators spans a dense linear subspace of $\ell _\infty (X)$, we see that $\ell _\infty (X)\subseteq \cstu (X)$.

Alternatively, the uniform Roe algebra of $X$ can be defined as the norm closure of the set of all bounded operators on $\ell _2(X)$ with \emph{finite propagation}, where the \emph{propagation} of an operator $a\in \cB (\ell _2(X))$ is given by
  \begin{equation}
  \label{EqDefinitionPropagartionOp}
  \propg (a)=\sup \{d(x,y)\ :\ \langle a\delta _y,\delta _x\rangle \not =0\}.
  \end{equation}

  \subsection{Flows on uniform Roe algebras}\label{SubsectionFlowDef}

With the present subsection we begin our detailed study of flows.

  \begin{definition}
  Given a $\mathrm C^*$-algebra $A$, by a \emph{pre-flow} on $A$, also known as a \emph{one-parameter group of automorphisms} of $A$, we shall mean a group homomorphism
  \[
  \sigma \colon {\mathbb R}\to \text {Aut(A)},
  \]
  from the  additive group of real numbers to the group of $^*$-automorphisms of $A$.  If the pre-flow $\sigma $ is \emph{strongly continuous} in the sense that $t\in \R \mapsto \sigma _t(a)\in A$ is continuous for all $a\in A$, then $\sigma $ is called a \emph{flow}.\footnote {Since the terminology ``pre-flow'' is not standard, we emphasize here that a pre-flow is simply a \emph{not necessarily strongly continuous flow}.}
  \end{definition}

We start our study of flows on uniform Roe algebras by showing that every such flow is completely determined by a self-adjoint operator on $\ell _2(X)$. Precisely, let $H$ be a Hilbert space and $h $ be a  self-adjoint operator on $H$.  Continuous functional calculus gives us that $e^{ih}$ is a unitary operator on $H$. Hence, each such operator induces a pre-flow $\sigma _h$ on $\cB (H)$ given by
  \begin{equation}
  \label{EqrefFlowIntro}
  \sigma _{h,t}(a)=e^{ith}ae^{-ith}
  \end{equation}
  for all $t\in \R $, and all $a\in \cB (H)$. It is well known that any weakly continuous flow on $\cB (H)$ is of this form for some self-adjoint operator $h$ on $H$ (see \cite [Example 3.2.35]{BratteliRobinsonVol1}).  In Subsection \ref{SubsectionFlows}, we employ this result to prove the following:

  \begin{proposition}
  \label{ThmFlowsGivenByH}
  Let $H$ be a separable Hilbert space, $ A \subseteq \cB (H)$ be a $\mathrm C^*$-algebra containing all compact operators, and $\sigma $ be a flow on $A$. Then there is a self-adjoint operator $h$ on $H$ such that
  \[
  \sigma _t(a)= e ^{ith}ae^{-ith}
  \]
  for all $t \in \R $, and all $a \in A$.
  \end{proposition}

Proposition \ref{ThmFlowsGivenByH} may be well known to experts, but for lack of a precise citing source and since this will play a very important role in our understanding of flows on uniform Roe algebras, we present a detailed proof in Subsection \ref{SubsectionFlows} below.

Since a uniform Roe algebra always contains the compacts, Proposition \ref{ThmFlowsGivenByH} applies to the flows that we are interested in.  Therefore it suffices to study flows on uniform Roe algebras arising from  self-adjoint operators on $\ell _2(X)$ via the formula~\eqref{EqrefFlowIntro}. This problem splits naturally into two parts:
  \begin{question} \label{QuestionsIandII}
  Let $h$ be a  self-adjoint operator on $\ell _2(X)$.
  \begin{enumerate}
  [label=\textnormal {(\Roman *)}]
    \item \label{Item2Intro} Under which conditions can we guarantee that $\sigma _{h,t} (\cstu (X))= \cstu (X)$ for all $t\in \R $?
    \item \label{Item1Intro} Assuming that the above holds, when is $\sigma _{h}$ strongly continuous on $\cstu (X)$? In other words, when is a pre-flow on $\cstu (X)$ of the form $\sigma _h$, an actual flow?
  \end{enumerate}
  \end{question}

To deal with  Question
  \ref{QuestionsIandII}\ref{Item1Intro}, we start with a detailed analysis of coarse maps in Section \ref{SectionFlowsGivenByMaps} which culminates in a noncommutative version of coarseness for self-adjoint operators on $\ell _2(X)$ and a notion of coarseness for flows. Precisely:

  \begin{definition}
  \label{DefinitionCoarseOp}
  Let $X$ be a u.l.f.\ metric space and let $h$ be a  closed operator on $\ell _2(X)$.\footnote{Recall, a \emph{closed operator} on a Hilbert space is a partially defined, not necessarily bounded, operator whose graph is closed.} We call $h$ a \emph{coarse operator} if
  the domain of such $h$ is \emph{invariant under partial translations} (see Definition \ref{DefiInvPartialTranslation}), meaning
  \[
  v_f(\dom (h))\subseteq \dom (h), 
  \]
  and the commutator\footnote {Throughout, $[a,b]$ denotes the commutator $ab-ba$.} of $h$ and $v_f$ satisfies
  \[
  \| [ h ,v_f] \| < \infty \ \text{ and }\ [h,v_f]\in \cstu(X)
  \]
  for all partial translations $f$ on $X$; in the latter, we identify $[h,v_f]$ with its unique bounded extension to  the whole $\ell_2(X)$.
  \end{definition}

We should remark that the natural domain of the commutator operator $[h,v_f]=hv_f-v_fh$ is $\dom (h)$, because $\dom (h)$ is invariant under $v_f$ by hypothesis. Also, the invariance of $\dom (h)$ under the $v_f$ easily implies that $c_{00}(X)$\footnote{Here $c_{00}(X)$ denotes the vector subspace of $\ell_2(X)$ of finitely supported vectors.} is contained in the domain of $h$, and in fact $c_{00}(X)$ will be shown to be a core for $h$, in the sense that the graph of the restriction of $h$ to $c_{00}(X)$ is dense in the graph of $h$.  In Theorem \ref{Thmc00IsAGreatCore!}, we will show even more, namely that every coarse operator is in fact \emph{admissible} in the sense of Definition \ref{DefineAdmissible}, below, a crucial tool for several of the arguments in our main proofs.

  \begin{definition}
  \label{DefinitionCoarseFlow}
  Let $X$ be a u.l.f.\ metric space and $\sigma $ be a flow on $\cstu (X)$. We call $\sigma $ a \emph{coarse flow} if
  \[
  t\in \R \mapsto \sigma _{h,t}(v_f)\in \cstu (X)
  \]
  is differentiable for all partial translations $f$ of $X$.
  \end{definition}

A few words are in place here to justify the definitions above. First, as will be seen in Proposition \ref{ProphIsCoarseIFF} below, a map $h\colon X\to \R $ is coarse if and only if its interpretation as a diagonal operator on $\ell _2(X)$ is a coarse operator. Moreover, any flow on a uniform Roe algebra induced by a map $h\colon X\to \R $ is automatically a coarse flow (see Proposition \ref{PropIfFlowIsGivenByFucntionFinPropOpAreInDomOfInfinitesimalGenerator}). The study of coarse flows and of coarse operators is then a natural generalization of the commutative diagonal scenario.

With the definitions above in hand, we obtain the following noncommutative version of Proposition \ref{Prop2.1BragaExel}; which deals with Question \ref{QuestionsIandII}\ref{Item1Intro} above.

  \begin{proposition}
  \label{ProphCoarseOpImpliesFLow}
  Let $X$ be a u.l.f.\ metric space and $\sigma $ be a pre-flow on $\cstu (X)$.  Then $\sigma $ is a coarse flow if and only if it is induced by a coarse operator $h$ on $\ell _2(X)$.
  \end{proposition}

While, for the case of diagonal flows, there is no reason to worry about Question \ref{QuestionsIandII}\ref{Item2Intro}, the treatment of this question is substantially more subtle in the noncommutative case. For starters, we must introduce the \emph{unbounded part of $\cstu (X)$} to the picture.

\begin{definition} \label{DefineAdmissible} Let $X$ be a u.l.f.\ metric space and $h$ be an   operator on $\ell _2(X)$. 
We call $h$ \emph{admissible} if it is closed, $c_{00}(X)\subseteq \dom (h)$, and
    \[
    h\xi =\sum _{y\in X}\langle \xi ,\delta _y\rangle h\delta _y
    \]
   for all $\xi \in \dom (h)$. We say that $h$ has \emph{finite propagation} if it is admissible and \[\sup \{d(x,y) \ :\ \langle h\delta _x,\delta _y\rangle \neq 0\}<\infty .\] \end{definition}

We emphasize here that, whenever we say that an operator has finite propagation, it is assumed that we are speaking of an admissible operator, or else the condition on matrix entries above might not mean much.  We note that this definition generalizes the usual concept of finite propagation for bounded operators to the unbounded setting.

We next introduce the \emph{unbounded part} of $\cstu (X)$, a concept which will, as the name suggests, involve unbounded operators and, crucially, the difference $h-h'$, of two unbounded operators $h$ and $h'$.  This is a notoriously touchy issue not least because $h$ and $h'$ may have very different domains.  For this reason, except for rare occasions, we shall only dare to speak of $h-h'$ when there is a common core to $h$ and $h'$.  Since every admissible operator admits $c_{00}(X)$ as a core, taking the difference of such operators will not involve any risks.

  \begin{definition}
  \label{DefinitionUnboundedPartofURA}
  Let $X$ be a u.l.f.\ metric space and $h$ be an admissible operator on $\ell _2(X)$. We say that $h$ is in the \emph{unbounded part of $\cstu (X)$} if, for all $\eps >0$, there is a finite propagation operator   $h'$ on $\ell _2(X)$ such that $h-h'$ (defined on its natural domain $\dom (h)\cap \dom (h')$) is a bounded operator with norm at most $\eps $. We denote the set of all such operators by $\mathrm {C}^*_{u,unb}(X)$.\footnote {In these notes, we consider $\mathrm {C}^*_{u,unb}(X)$ simply as a set and do not attempt to view it with any algebraic or topological structure.}
  \end{definition}

Notice that any  operator on $\ell _2(X)$ with finite propagation is in $\mathrm {C}^*_{u,unb}(X)$. In particular, any diagonal operator, regardless of being given by a coarse map, is in $\mathrm {C}^*_{u,unb}(X)$. On the other hand, the intersection of $\cB (\ell _2(X))$ with $\mathrm {C}^*_{u,unb}(X)$ is precisely $\cstu (X)$.

  Another way of thinking of $\mathrm {C}^*_{u,unb}(X)$ is to consider the set $\Adm $ formed by all admissible operators, equipped with the (possibly infinite) metric
  \[
  d(h,k) = \|(h-k)\restriction _{c_{00}(X)}\|.
  \]
  Within it one has the subset $\mathrm {C}^*_{u,unb}[X]$ formed by the operators of finite propagation, and then
  $\mathrm {C}^*_{u,unb}(X)$
  turns out to be precisely the closure of
  $\mathrm {C}^*_{u,unb}[X]$
  in $\Adm $.

With the purpose of keeping this introduction to a minimum, we will now restrict our results to metric spaces with G. Yu's property A: this is a fairly general geometric property imposed on the u.l.f.\ metric space $X$ which is often seen as a sort of ``geometric amenability'' (see Subsection \ref{SubsectionCoarseIntro} for a precise definition), which is equivalent to nuclearity of $\cstu (X)$. For more general results which do not depend on property A, we simply refer the reader to the appropriate statements in the paper. For the specialists, we would like to say that the use of property A is needed so that the uniform Roe algebra of $X$ equals its \emph{quasi-local algebra} (see Subsection \ref{SubsectionCoarseIntro} for definitions).

For u.l.f.\ metric spaces with property A, we completely characterize the self-adjoint operators which give rise to coarse flows on uniform Roe algebras. Precisely, we obtain the following.

  \begin{theorem}
  \label{ThmCharactCoarseFlowsInTermsOp}
  Let $X$ be a u.l.f.\ metric space with property A, $h$ be a  self-adjoint operator on $\ell _2(X)$, and $\sigma _h$ be the pre-flow on $\cB (\ell _2(X))$ induced by $h$. The following are equivalent.
  \begin{enumerate}
         \item \label{ThmCharactCoarseFlowsInTermsOpItem2} $\sigma _h$ is a coarse flow on $\cstu (X)$.
         \item \label{ThmCharactCoarseFlowsInTermsOpItem1} $h$ is a coarse operator and it belongs to $\mathrm {C}^*_{u,unb}(X)$.
         \item \label{ThmCharactCoarseFlowsInTermsOpItem3} The map $E(h)\colon X\to \R $ given by
  \[
  E(h)_x=\langle h\delta _x,\delta _x\rangle\ \text { for all } \ x\in X
  \]
  is coarse and, interpreting $E(h)$ as a diagonal operator, $h-E(h)$ is a bounded operator in $\cstu (X)$.\footnote {The notation ``$E(h)$'' is an allusion to the canonical conditional expectation from $\cB (\ell _2(X))\to \ell _\infty (X)$, the formula of which is given analogously.}
  \end{enumerate}
  \end{theorem}

The main ingredient for the implication \eqref{ThmCharactCoarseFlowsInTermsOpItem1}$\Rightarrow $\eqref{ThmCharactCoarseFlowsInTermsOpItem2} is Theorem \ref{ThmCoarseOpInUmbPart}, which shows that for any such $h$, the bounded operator $e^{ih}$ is in the quasi-local algebra of $X$ for all $t\in {\mathbb R}$. On the other hand, the implication \eqref{ThmCharactCoarseFlowsInTermsOpItem2}$\Rightarrow $\eqref{ThmCharactCoarseFlowsInTermsOpItem1} of Theorem \ref{ThmCharactCoarseFlowsInTermsOp} is always valid regardless of property A. This is obtained as Theorem \ref{ThmhCoarseIsInUnboundedPartOfURA} below and its proof uses an adaptation to the unbounded setting of the machinery developed in \cite {LorentzWillett2020} to show that bounded derivations on uniform Roe algebras are always inner.

  \subsection{Cocycle equivalence of flows}

Once Items \ref{Item2Intro} and \ref{Item1Intro} of Question \ref{QuestionsIandII} are dealt with, we move to classifying when coarse flows are, in some sense, equivalent:

\begin{question}
Let $h$ and $k$ be self-adjoint operators on $\ell _2(X)$   inducing flows $\sigma_h$ and $\sigma_k$ on $\cstu(X)$.
  \begin{enumerate}  
  [label=\textnormal {(\Roman *)}]\setcounter {enumi}{2}
    \item What relation must be satisfied between $h$ and $k$ so that the flows $\sigma _h$ and $ \sigma _k$ are ``morally'' the same?
  \end{enumerate}
  \end{question}
This will be formalized in terms of \emph{cocycle equivalence} between flows.
For this, the notion of cocycle will be essential. Recall that we have the following well-known equivalence relations among flows on $\mathrm C^*$-algebras:

  \begin{definition}
  \label{DefiCocycleConjugate}
  Let $A$ be a unital $\mathrm C^*$-algebra, $\sigma $ be a pre-flow on $A$, and $(u_{t})_{t\in \R }$ be a norm continuous family of unitary operators in $A$.
  \begin{enumerate}
  \item We say that $(u_{t})_{t\in \R }$ is a \emph{cocycle for
  $\sigma $} if
  \[
  u_{t+s}=u_t\sigma _t(u_s) \text { for all }\ t,s\in \R .
  \]
  \item If $(u_t)_{t\in \R }$ is a cocycle for $\sigma $, then the expression
  $\rho _t(a)=u_t\sigma _t(a)u^*_t$, for $a\in A$ and $t\in \R $, defines a pre-flow $\rho $ on $A$. We call $\rho $ a \emph{cocycle perturbation of $\sigma $}.
  \item If, moreover, there exists a family $(u_t)_{t\in \R }$, as above, that is differentiable at $0$, then $\rho $ is called an \emph{inner perturbation of $\sigma $}.\footnote {If this is the case, $(u_t)_{t\in \R }$ is actually differentiable everywhere, and we refer the reader to \cite [Section 1]{Kishimoto2000RepMathPhy} for further details.}
  \item We say that $\sigma $ is \emph{cocycle conjugate to $\rho $} if there is an automorphism $\tau $ of $A$ such that $\rho $ is a cocycle perturbation of $\tau \circ \sigma \circ \tau ^{-1}$, where the later is the pre-flow given by $t\in \R \mapsto \tau \circ \sigma _t\circ \tau ^{-1}\in \mathrm {Aut}(A)$.
  \end{enumerate}
  \end{definition}

All these notions define equivalence relations on the set of pre-flows on a given $\mathrm C^*$-algebra and relate to each other as in the diagram below. Moreover, if a pre-flow is equivalent under any of these notions to an actual flow, then it must be a flow as well.\\

  \begin{tikzpicture}
  \xymatrix {\mathrm {Inner\: perturbation} \ar [r] & \mathrm {Cocycle\: perturbation} \ar [r] & \mathrm {Cocycle\: conjugacy}. }
  \end{tikzpicture}\\

For metric spaces with property A, we are capable of showing that every coarse flow is, up to cocycle perturbation, a diagonal flow. Moreover, this is more than an ``existence theorem'' as the diagonal cocycle perturbation can be actually computed in a fairly natural way.  As diagonal flows are much easier to deal with, this can drastically reduce any analysis needed when working with coarse flows. In particular, this allows the machinery developed in \cite {BragaExel2023KMS} for the study of KMS states of diagonal flows to be applied to more general flows. Precisely, we prove the following.

  \begin{theorem}
  \label{ThmFinPropagationCocycleEquivDiagFlowPropA}
  Let $X$ be a u.l.f.\ metric space with
  property A, and let $h$ be a coarse
operator on $\ell_2(X)$ such that
  \[
  e^{ith}\cstu (X)e^{-ith}=\cstu (X)
  \]
  for all $t\in {\mathbb R}$.
  Then the flow $\sigma _h$  on $\cstu (X)$ induced by $h$ is a
cocycle perturbation of a diagonal flow.  Indeed, $\sigma _h $ is a cocycle perturbation of the flow $\sigma _{E(h)}$, where $E(h)\colon X\to \R $ is the coarse map given by
  \[
  x\in X\mapsto \langle h\delta _x,\delta _x\rangle \in \R\ \text { for all }\ x\in X.
  \]
  \end{theorem}

Besides reducing the problem to the diagonal case, we also characterize when self-adjoint operators on $\ell _2(X)$ give rise to ``equivalent'' flows. For the commutative case, no geometric condition on $X$ is needed, and we obtain not only a classification of flows, but also of pre-flows (notice here that the property of being a flow is preserved even under the weakest form of equivalence given by cocycle conjugacy):

  \begin{theorem}
  \label{ThmCocycleEquivFlowCoarseMaps}
  Let $(X,d)$ be a u.l.f.\ metric space, and let $h$ and $k$ be real-valued maps on $X$. The following are equivalent:
  \begin{enumerate}
  \item \label{ThmCocycleEquivFlowCoarseMapsItem1} The maps $h$ and $k$ are close, i.e., $\sup _{x\in X} d(h_x,k_x)<\infty $,
  \item \label{ThmCocycleEquivFlowCoarseMapsItem2} The pre-flows $\sigma _h$ and $\sigma _k$ are inner perturbations of each other.
  \item \label{ThmCocycleEquivFlowCoarseMapsItem3} The pre-flows $\sigma _h$ and $\sigma _k$ are cocycle perturbations of each other.
  \end{enumerate}
  \end{theorem}

Notice that the property in Theorem \ref{ThmCocycleEquivFlowCoarseMaps}\eqref{ThmCocycleEquivFlowCoarseMapsItem1} above is equivalent to the difference operator $h-k$, given by the interpretation of the maps $h,k\colon X\to \R $ as operators on $\ell _2(X)$, being a bounded operator.

In order to obtain Theorems \ref{ThmFinPropagationCocycleEquivDiagFlowPropA} and \ref{ThmCocycleEquivFlowCoarseMaps}, we need yet another important ingredient which is interesting on its own: 

\begin{theorem} \label{ThmCoroFlowsuRaSpaciallyImplemented} Let $X$ be a u.l.f.\ metric space with property A, and let $h$ be a   self-adjoint operator on $\ell_2(X)$ such that $\sigma _h$ induces a flow on $\cstu (X)$.  Then $e^{ith}\in \cstu (X)$ for all $t\in {\mathbb R}$.  \end{theorem}

In fact, as the reader will see in Section \ref{SectionAutConnectedIdentity}, this is not a result about flows, but merely about strongly continuous paths of automorphisms. The proof of this result is quite technical and it makes use of three ingredients: (1) basics of K-theory, (2) techniques developed in \cite {WhiteWillett2017} for the study of Cartan subalgebras of uniform Roe algebras, and (3) techniques developed for the solution of the rigidity problem of uniform Roe algebras (\cite {BragaFarah2018Trans,BaudierBragaFarahKhukhroVignatiWillett2021uRaRig}).

For the general noncommutative case, we must  once again restrict ourselves to metric spaces
with property A. Our main results are the following.

  \begin{theorem}
  \label{ThmClassificationFlowsGenMaisFraco} Let $X$ be a u.l.f.\ metric space with property A, and let $h$ and $k$ be coarse operators on $\ell _2(X)$ inducing flows $\sigma _h$ and $\sigma _k$ on $\cstu (X)$. If   $h-k$ is bounded on $c_{00}(X)$,
  then the flows $\sigma _h$ and $\sigma _k$ on $\cstu (X)$ are cocycle perturbations of each other.
  \end{theorem}

  \begin{theorem}
  \label{ThmClassificationFlowsGen} Let $X$ be a u.l.f.\ metric space with property A, and let $h$ and $k$ be coarse,  self-adjoint operators on $\ell _2(X)$, whose induced flows on $\cstu (X)$ are respectively denoted by $\sigma _h$ and $\sigma _k$. The following are equivalent:
  \begin{enumerate}
    \item \label{ThmClassificationFlowsGenItem1} There is an unitary $u\in \cstu
(X)$ such that $\dom (h)=\dom (uku^*)$, and
$h-uku^*$ is bounded on this common domain.
    \item \label{ThmClassificationFlowsGenItem2} The flows $\sigma _h$ and $\sigma _k$ are cocycle conjugate to each other.
  \end{enumerate}
  \end{theorem}

Finally, in Section \ref{SectionPreflows}, we show that, without the condition of property A, and working only with pre-flows, the results above do not hold. More precisely, we show that in any expander graph, one can find a natural pre-flow on its uniform Roe algebra which is not a cocycle perturbation of a diagonal flow (see Theorem \ref{ThmFlowExpGraph}).

\section{Preliminaries}\setcounter{subsection}{1}

In this section, we recall the basic tools that we will need about unbounded operators, one-parameter unitary groups, flows on algebras of operators, and coarse geometry.  For a detailed treatment of these, we refer the reader to \cite [Chapter 5]{Pede:Analysis},  \cite [Chapter 13]{rudin1991functional}  or \cite{BratteliRobinsonVol1} for the first topics, and to \cite {RoeBook,NowakYuBook} for the last one. Besides setting up terminology, we also refine a well-known result on weakly continuous flows on $\cB (H)$ \cite [Example 3.2.35]{BratteliRobinsonVol1} to guarantee that any flow on a uniform Roe algebra is determined by an unbounded operator (Proposition \ref{ThmFlowsGivenByH}).

  \subsection{ Unbounded Operators on Hilbert spaces}

 Let $H$ be a Hilbert space. From now on, the term \emph{operator on $H$} will refer to a (not necessarily bounded) linear function
  \[
  h\colon \dom (h)\subseteq H\to H,
  \]
  where $\dom (h)$ is a vector subspace of $H$. The unbounded operators of interest
in this paper will be in general self-adjoint. In particular, these operators are
automatically closed and densely defined, i.e., $\dom (h)$ is dense in $H$. Addition,  subtraction and multiplication of operators are defined in the obvious ways taking the necessary precautions regarding domains; precisely,
  \[
  \dom (h\pm k)=\dom (h)\cap \dom (k)
  \]
  and
  \[
  \dom (kh)=\{\xi \in \dom (h)\ :\ h\xi \in \dom (k)\}.
  \]
   The space of all (globally defined) bounded operators on $H$ is denoted by $\cB (H)$ and its ideal of compact operators by $\cK (H)$.

Given the fundamental role played by the matrix representation of operators on $\ell_2(X)$ in the study of $\cstu (X)$, e.g., in the definiton of the notion of propagation, it is natural that unbounded operators also be seen from the point of view of matrices.  However one cannot expect to make sense out of the matrix entries
  \[
  h_{x,y} = \langle h\delta _y,\delta _x\rangle
  \]
  of an unbounded operator $h$ on $\ell_2(X)$, since such operators are only densely defined, and there is no reason for $\delta _y$ to lie in its domain.  In fact, even if all basis vectors $\delta _y$ belong to $\dom (h)$, the values of $h$ on these vectors might not necessarily determine $h$ outside $c_{00}(X)$.  After all, the fact that $h$ is discontinous, as a rule, renders the density of $c_{00}(X)$ within $\ell_2(X)$ useless for the purpose of inferring the value of $h(\xi )$, when $\xi $ is not in $c_{00}(X)$.

The operators that will play important roles in this work will therefore all be \emph{admissible} in the following sense:

\begin{definition} [cf.\ Definition \ref{DefineAdmissible}]
  Let $X$ be any set.
  An operator $h$ on $\ell_2(X)$ will be called \emph{admissible} if it is closed, $c_{00}(X)\subseteq \dom (h)$, and
    \[
    h\xi =\sum _{y\in X}\langle \xi ,\delta _y\rangle h\delta _y
    \]
   for all $\xi \in \dom (h)$.  \end{definition}

Under the above conditions, the coordinates of $h\xi $ relative to each basis vector $\delta _x$ are then given by
  \[
  \langle h\xi , \delta _x\rangle =
  \sum _{y\in X}\langle \xi ,\delta _y\rangle \langle h\delta _y,\delta _x\rangle =
  \sum _{y\in X}h_{x, y}\xi _y,
  \]
  where $h_{x, y} = \langle h\delta _y, \delta _x\rangle $, and $\xi _y=\langle \xi ,\delta _y\rangle $.
  In other words, the matrix representation of $h$ holds, not only for vectors in $c_{00}(X)$, but also for every vector $\xi $ in the domain of definition of $h$. Furthermore, given any pair $(\xi , h\xi )$ in the graph of $h$, we may write
  $ \xi =\lim _F \xi _F$,
  where
  \[
  \xi _F=\sum _{y\in F}\langle \xi , \delta _y\rangle \delta _y,
  \]
  the limit taken for $F$ running in the directed set of all finite subsets $F\subseteq X$, and then the fact that $h$ is admissible  implies that $\big (\xi _F,h(\xi _F)\big )$ converges to $(h\xi ,\xi )$, which is to say that $c_{00}(X)$ is a core for $h$.

   Admissibility is a very useful property since it implies that $c_{00}(X)$ is not only a core, but it is actually a  \emph{`great core'} in the sense that we can choose our approximating sequences in the most natural way. Precisely, firstly recall that a subspace $D\subseteq \dom (h)$ is called a \emph{core for $h$} if the graph of $h\restriction D$ is dense in the graph of $h$. In other words, if $\xi \in \dom (h )$, then there is a sequence $(\xi _n)_n$ in $D$, such that $\xi =\lim _n\xi _n$, and $h\xi =\lim _nh\xi _n$. Admissibility says that the sequence $(\xi _n)_n$ can be simply taken to be the truncations of $\xi $ with respect to the canonical basis of $\ell _2(X)$.

We finish this subsection with a useful elementary result.

    \begin{proposition} \label{Prop05Nov24}
Let $X$ be a u.l.f.\ metric space and let  $h$ and $k$ be admissible operators on $\ell_2(X)$.  The following holds.
  \begin{enumerate}
  [label=\textnormal {(\alph *)}]
  \item
  $
  \vrule height 12pt depth 8pt width 0pt 
  \|h-k\| = \|(h-k){\restriction}_{c_{00}(X)}\|$.\footnote{Here $\|h\|:=\sup\{\|h(\xi)\|: \|\xi\|\leq 1, \xi\in \dom(h)\}$ denotes the possibly infinite norm of the (unbounded) operator $h$.} 
  \item If $h-k$ is bounded on $c_{00}(X)$, then $\dom (h)= \dom (k)$.
  \end{enumerate}
  \end{proposition}

\begin{proof}
  The first point is an easy consequence of the fact that $c_{00}(X)$ is a common core for $h$ and $k$, so we leave the details for the reader.

As for the second point, let $\xi\in \dom(h)$ and let us show that $\xi\in\dom(k)$. For this choose a increasing sequence $(F_n)_n$ of finite subsets of $X$, whose union coincides with $X$, and for each $n$, let $\xi_n$ be the projection of $\xi$ on $\ell _2(F_n)$.  Then clearly $\xi_n\to \xi $ and, since $h$ is admissible, also $h\xi_n\to h\xi$.  Noticing that $(\xi_n)_n$ is a sequence in $c_{00}(X)$, and that $h-k$ is bounded on this subspace,  we have that $((h-k)\xi _n)_n$ is a Cauchy sequence, hence so is $(k\xi _n)_n$.  As $k$ is closed, this shows that $\xi\in \dom(k)$. This proves $\dom(h)\subseteq \dom(k)$ and, by symmetry of the hypothesis, the inclusion $\dom(k)\subseteq \dom(h)$ follows analogously.  \end{proof}

  \subsection{Flows and one-parameter groups}\label{SubsectionFlows}

  \begin{definition}
   Let $H$ be a Hilbert space and $(u_t)_{t\in \R }$ be a family of unitary operators in $\cB (H)$. We call $(u_t)_{t\in \R }$ a \emph{one-parameter unitary group on $H$} if
  \[
  u_{t+s}=u_tu_s\ \text { for all }\ t,s\in \R .
  \]
  If, moreover,
  \[
  t\in \R \mapsto u_t\xi \in H
  \]
  is continuous for all $\xi \in H$, we say that
  $(u_t)_{t\in \R }$ is \emph{strongly continuous}.
  \end{definition}

Given a self-adjoint operator $h$ on a Hilbert space $H$, the spectral theorem allows us to define $e^{ih}$, a bounded operator on $H$. Moreover,
  \[
  t\in \R \mapsto e^{ith}\in \cB (H)
  \]
  is a strongly continuous one-parameter unitary group and
  \begin{equation}
  \label{EqDifOneParGroupGivenSAOp}
  \lim _{t\to 0}\frac {e^{ith}\xi -\xi }{t}=ih\xi
  \end{equation}
  for all $\xi \in \dom (h)$.
  Conversely,  if $\xi \in H$ is such that the limit in \eqref{EqDifOneParGroupGivenSAOp} exists,
then necessarily $\xi \in \dom (h)$, and the identity in \eqref{EqDifOneParGroupGivenSAOp}
holds
  (see \cite [Proposition 5.3.13]{Pede:Analysis}).  Moreover, by Stone's theorem (see \cite [Theorem 5.3.15]{Pede:Analysis}), every strongly continuous one-para\-meter unitary group $t\in \R \mapsto u_t\in \cB (H)$ is of the above form, that is, $u_t=e^{ith}$ for all $t\in \R $, for a unique self-adjoint operator $h$ on $H$.

We now show that any flow on a uniform Roe algebra is given by conjugating by a one-parameter unitary group; this is a refinement of \cite [Example 3.2.35]{BratteliRobinsonVol1} for our setting. This will be essential in our study of flows. Parts of the proof of Proposition \ref{ThmFlowsGivenByH} below are well known to experts, but we provide details for clarity and for the non-experts' convinience.

  \begin{proof}
  [Proof of Proposition \ref{ThmFlowsGivenByH}] We start noticing that $\cK (H)$ is invariant under $\sigma _t$. For that, observe that $\cK (H)$ is the unique minimal ideal of $A$. Indeed, if $J$ is another nonzero ideal of $A$, we have
  \[
  J \cap \cK (H) = J \cdot \cK (H) \neq \{0\},
  \]
  since $\cK (H)$ is an essential ideal. As
  $\cK (H)$ is a simple $\mathrm C^*$-algebra, it follows that $J \cap \cK (H) = \cK (H)$, whence $\cK (H) \subseteq J$.  As $\cK (H)$ admits the above characterization in precise algebraic terms, we conclude that $\cK (H)$ is invariant under automorphisms of $A$; in particular, $\sigma _t(\cK (H)) = \cK (H)$ for every $t \in \R $.

For each $t\in \R $, let $\tau _t$ be the restriction of $\sigma _t$ to $\cK (H)$, so that $\tau $ becomes a flow on $\cK (H)$.  Observing that $\B (H)$ is the multiplier algebra of $\cK (H)$, each $\tau _t$ extends canonically to an automorphism $\rho _t$ of $\cB (H)$, and it is not difficult to show that $\rho _t\rho _s = \rho _{t+s}$ for all $t, s\in \R $, i.e., that
  $\rho $ is a pre-flow on $\cB (H)$. Observe that, for $t\in \R $, $a\in A$, and $k\in \cK (H)$, we have
  \[
  \rho _t(a)k = \rho _t(a\rho _{-t}(k)) = \tau _t(a\tau _{-t}(k)) = \sigma _t(a\sigma _{-t}(k)) = \sigma _t(a)k.
  \]
  So, $\rho _t(a) = \sigma _t(a)$. Hence, $\rho $ extends $\sigma $ from $A$ to $\cB (H)$.
  Let us now show that $\rho $ is a weakly continuous flow on $\cB (H)$.  For this, pick any vector $\xi $ in $H$, and let $p_\xi $ be the projection onto the one-dimensional space spanned by $\xi $.  For every bounded operator $b$ on $H$ we then have that 
  \[
  \rho _t(b)\xi =
  \rho _t(b)p_\xi \xi =
  \rho _t\big (b\rho _{-t}(p_\xi )\big )\xi =
  \sigma _t\big (b\sigma _{-t}(p_\xi )\big )\xi ,
  \]
  because both $p_\xi $ and $b\rho _{-t}(p_\xi )$ are compact operators, and hence belong to $A$.  Using that $\sigma $ is a (continuous) flow on $A$, we see that $\rho _t(b)\xi $ is continuous as a function of $t$, meaning that $\rho $ is strongly continuous, and hence also weakly continuous.

Now, \cite [Example 3.2.35]{BratteliRobinsonVol1} shows that a weakly
continuous flow on $\cB (H)$ must be of the desired form,  so the result follows.
  \end{proof}

  \subsection{Coarse geometry}\label{SubsectionCoarseIntro}

Although it has already been briefly presented in the introduction, we now properly define the morphisms and equivalences in the coarse category.  Given metric spaces $(X,d_X)$ and $(Y,d_Y)$, a map $f\colon X\to Y$ is called \emph{coarse} if for all $r>0$,
  \[
  \sup _{d_X(x,z)\leq r} d_Y(f(x),f(z))<\infty .\]
  Maps $f,g\colon X\to Y$ are called \emph{close} if
  \[
  \sup _{x\in X}d_Y(f(x),g(x))<\infty .
  \]
  If $f\colon X\to Y$ is coarse and there is a coarse $g\colon Y\to X$ such that $g\circ f$ and $f\circ g$ are close to $\mathrm {id}_{X}$ and $\mathrm {id}_Y$, respectively, then $f$ is called a
  \emph{coarse equivalence}. 

The uniform Roe algebra has a very natural dense $^*$-subalgebra which will be useful for us. Precisely, the algebraic counterpart of $\cstu (X)$, denoted $\cstu [X]$, is the $^*$-subalgebra of $\cB (\ell _2(X))$ generated by all $av_f$, where $a\in \ell _\infty (X)$, and $f$ is a partial translation of $X$. Unless $X$ is finite, $\cstu [X]$ will never be closed.

As already made clear in the introduction, the metric spaces of interest to us will be \emph{uniformly locally finite}, or \emph{u.l.f.} for short (see Subsection \ref{SubsectionURA}
  for definition). One of the great advantages of such metric spaces is that an $X$-by-$X$ matrix of complex numbers with finite propagation induces a bounded operator on $\ell _2(X)$ if and only if its entries are uniformly bounded. For further reference, we isolate this as a proposition below (see \cite [Lemma 4.27]{RoeBook}).

  \begin{proposition}
  \label{PropOperatorInducedFinPropULF}
  Let $X$ be a u.l.f.\ metric space. Then, for all $r>0$, if $A=[a_{x,y}]_{x,y\in X}$ is an $X$-by-$X$ matrix of complex numbers with
  \[
  \sup _{x,y\in X} |a_{xy}|<\infty,
  \]
  and $\propg (A)\leq r$, then $A$ canonically induces a bounded operator $a$ on $\ell _2(X)$, such that $\langle a\delta _y,\delta _x\rangle =a_{xy}$ for all $x,y\in X$. Moreover, there is a positive constant $M$, dependig only on $X$ and $r$, such that
  \[
  \|a\|\leq M \sup _{x,y\in X} |a_{xy}|.
  \]
  \end{proposition}

An important class of metric spaces in coarse geometry is the one of u.l.f.\ metric spaces with property A, introduced by G. Yu in \cite {Yu2000}. Although the definition of property A will not be used per se, we quickly recall it here since this property will play a major role throughout these notes.

  \begin{definition}
    Let $(X,d)$ be a u.l.f.\ metric space. We say that $X$ has \emph{property A} if, for all $\eps ,r>0$, there is a function $\xi \colon X\to \partial B_{\ell _2(X)}$\footnote {Here $\partial B_{\ell _2(X)}$ denotes the unit sphere of $\ell _2(X)$.} such that
  \begin{enumerate}
        \item $\|\xi _x-\xi _y\|\leq \eps $ for all $x,y\in X$ with $d(x,y)\leq r$, and
        \item $\sup \{d(x,y)\ :\ \langle \xi _x,\delta _y\rangle \neq 0\}<\infty $.
    \end{enumerate}
  \end{definition}

  It is well known that a u.l.f.\ metric space $X$ has property A if and only if $\cstu (X)$ is a nuclear $\mathrm C^*$-algebra (see \cite [Theorem 5.5.7]{BrownOzawa}).

  Besides the uniform Roe algebra of $X$, the quasi-local algebra of $X$ will play an important role in these notes.

  \begin{definition}
     Let $(X,d)$ be a u.l.f.\ metric space. The \emph{quasi-local algebra of $X$}, denoted by $\cstql (X)$, is the $\mathrm C^*$-subalgebra of $\cB (\ell _2(X))$ consisting of all operators $a\in \cB (\ell _2(X))$ such that, for all $\eps >0$, there is $r>0$, such that for all $A,B\subseteq X$,
  \[
  d(A,B)>r\ \text { implies }\
  \|p_Aap_B\|\leq \eps .
  \footnote { Recall, as stated in the introduction, that $p_A$ is the orthogonal projection from $\ell _2(X)$ to $\ell _2(A)$.}
  \]
  Equivalently, a bounded operator $a$ lies in $\cstql (X)$, if and only if
  \[
  \lim _{r\to \infty }\ \sup _{\scriptscriptstyle d(\kern -1pt A,B)>r} \|p_Aap_B\| = 0.
  \]
  \end{definition}

While it is immediate that $\cstu (X)\subseteq \cstql (X)$ for any u.l.f.\ metric space $X$, it was an important open problem whether $\cstu (X)$ and $\cstql (X)$ coincide. J.\ \v {S}pakula and J.\ Zhang proved that this is indeed the case for u.l.f.\ metric spaces with property A (see \cite [Theorem 3.3]{SpakulaZhang2020JFA}). Very recently, N.\ Ozawa has shown that these algebras are never the same as long as $X$ contains a coarse copy of an expander graph (see \cite [Corollary C]{Ozawa2023uRaSmallerQL}).

\section{Diagonal flows}
  \label{SectionFlowsGivenByMaps}\setcounter{subsection}{1}

The goal of this section is to set the ground needed to define a noncommutative version of coarseness in the next section. More precisely, we characterize coarseness of a real valued function on a metric space $X$ in terms which allow us to define ``coarseness'' for an arbitrary self-adjoint operator on $\ell _2(X)$ (Proposition \ref{ProphIsCoarseIFF}). We also pin-down a property satisfied by flows induced by maps $X\to \R $ which will be needed in our study of more general flows induced by self-adjoint operators on $\ell _2(X)$ (Proposition \ref{PropIfFlowIsGivenByFucntionFinPropOpAreInDomOfInfinitesimalGenerator}).

  Our motivation for wanting to characterize coarseness of maps $X\to \R $ comes from the fact that coarse maps are precisely the ones which induce flows on $\cstu (X)$ in a canonical way; which in turn was intensively studied in \cite {BragaExel2023KMS}. Let us start by explaining this carefully.  An arbitrary function $h\colon X\to \R $ canonically induces a self-adjoint operator on $\ell _2(X)$ by multiplication. By abuse of notation, we denote this operator by $h$ as well. Precisely, $h$ is the operator with
  \[
  \dom (h)=\left \{\xi \in \ell _2(X)\ :\ \sum _{x\in X}|h_x\xi _x|^2<\infty \right \},
  \]
  and given by the formula
  \[
  h\xi =(h_x\xi _x)_{x\in X}
  \]
  for all $\xi \in \dom (h)$. Letting $c_{00}(X)$ denote the subspace of $\ell _2(X)$ formed by the vectors with finitely many nonzero coordinates, it is clear that
  $c_{00}(X)\subseteq \dom (h)$, so $\dom (h) $ is dense in $\ell _2(X)$. See Proposition \ref{PropCoarseImpliesDomPreervation} below for a proof that $h$ is actually an admissible operator. Moreover, $h$ defines a bounded operator on $\ell _2(X)$ if and only if $h\colon X\to \R $ is bounded. Hence, the identification above of functions $X\to \R $ with possibly unbounded self-adjoint operators on $\ell _2(X)$ is compatible with the usual identification of $\ell _\infty (X)$ as the bounded multiplication operators on $\ell _2(X)$, i.e., the diagonal operators with respect to the canonical orthonormal basis $(\delta _x)_{x\in X}$ of $\ell _2(X)$.

  \begin{convention}
  In order to avoid any confusion  stemming from the abuse of notation of denoting by the same symbol --- $h$ usually --- both a map $ X\to \R $ and its interpretation as an unbounded operator on $\ell _2(X)$, we write $h=(h_x)_{x\in X}$ when $h$ is being seen as a function $h\colon x\in X\mapsto h_x\in \R $. If $h$ is being interpreted as an operator on $\ell _2(X)$, its evaluation at a vector $\xi \in \ell _2(X)$ is written as $h\xi $. Also, $\dom (h)$ will always refer to the domain of $h$ as an operator on $\ell _2(X)$ (the domain of $h$ as a function $X\to \R $ is simply $X$ and this will not be a point of relevance in these notes).
  \end{convention}

Suppose that $(X,d)$ is a u.l.f.\ metric space and that $h\colon X\to \R $ is a map. Then, $e^{ih}$ may be seen as a bounded operator on $\ell _2(X)$ with propagation zero, so, $e^{ih}\in \cstu (X)$. In particular, $e^{ith}ae^{-ith}\in \cstu (X)$ for all $t\in \R $, and all $a\in \cstu (X)$, so the map
  \[
  \sigma _{h, t}\colon a\in \cstu (X) \mapsto e^{ith}ae^{-ith}\in \cstu (X)
  \]
  is an automorphism of $\cstu (X)$, and it can be easily shown to define a
pre-flow $\sigma _h$ on $\cstu (X)$. As seen in Proposition \ref{Prop2.1BragaExel},
$\sigma _h$ is a flow if and only if the map $h\colon X\to \R $ is coarse,  in which case we call $\sigma _h$ a \emph{diagonal flow}.

Proposition \ref{Prop2.1BragaExel}
  is one of our hints to finding an appropriate operator version of coarseness.
  Before that, we must notice that the domain of an operator induced by a coarse map $X\to \R $ is invariant with respect to the coarse geometry of $X$. We introduce the following definition.

  \begin{definition}
  \label{DefiInvPartialTranslation} Let $E\subseteq \ell _2(X)$ be a linear subspace.  Given a partial translation $f$ of $X$, we say that $E$ is \emph{$f$-invariant} if $v_f(E)\subseteq E$. If $E$ is $f$-invariant for all partial translations $f$ of $X$, we say that $E$ is \emph{invariant under partial translations}.
  \end{definition}

  \begin{remark}
  \label{RemarkDomainhInvariantContainsc00}
  Notice that if the domain of a self-adjoint operator $h$ on $\ell _2(X)$ is invariant under partial translations, then $\dom (h)$ automatically contains $c_{00}(X)$. Indeed, as $h$ is self-adjoint, its domain is densely defined and, in particular,   non-zero. Pick a non-zero $\xi\in \mathrm{dom}(h)$  and $x_0\in X$ such that the $x_0$-coordinate of $\xi$ is non-zero. Then, given $x\in X$ and letting  $f_{x}\colon\{x_0\}\to \{x\}$ be the partial translation taking $x_0$ to $x$, we have that    $v_{f_x}\xi$ is in $\mathrm{dom}(h)$ for all $x\in X$. As $v_{f_x}\xi$ is a non-zero vector supported on $x$ and $x\in X$ was arbitrary, we conclude that $\mathrm{dom}(h)$ contains $c_{00}(X)$. This will be used throughout these notes without further mention. The fact that $c_{00}(X)\subseteq \dom (h)$ does not necessarily mean that $h$ is admissible, although most of the operators in this work will possess this property.
  \end{remark}

  For simplicity, we introduce the following (almost) self-explaining terminology: given a partial translation $f$ of $X$, and $r\geq 0$, we say that $f$ is a \emph{partial $r$-translation} if
  \[
  \sup _{x\in \dom (f)}d(f(x),x)\leq r.
  \]

  \begin{proposition}
  \label{PropCoarseImpliesDomPreervation}
  Let $X$ be a u.l.f.\ metric space and let $h\colon X\to \R $ be any function. Then, seen as a multiplication operator, $h$ is admissible.  Moreover, if $h$ is a coarse function, then $\dom (h)$, the domain of $h$ as an operator on $\ell _2(X)$, is invariant under partial translations.
  \end{proposition}

  \begin{proof}It is evident that $\dom (h)$ contains $c_{00}(X)$, and we leave it for the reader to prove the standard fact that $h$ is a closed operator.  On the other hand, given $\xi \in \dom (h)$, we have that
  \[
  h\xi =(h_x\xi _x)_{x\in X} = \sum _{x\in X}h_x\xi _x\delta _x = \sum _{x\in X}\langle \xi , \delta _x\rangle h\delta _x,
  \]
  so $h$ is admissible.  In order to show that $\dom (h)$ is invariant under partial translations, suppose that
  $h\colon X\to \R $ is coarse and that $\xi =(\xi _x)_{x\in X}$ is in the domain of $h$, seen as an operator on $\ell _2(X)$. So,
  \[
  \| h\xi \|^2=\sum _{x\in X}|h_x\xi _x|^2<\infty .
  \] Let $f$ be a partial translation of $X$. As $h$ is coarse,
  \[
  s := \sup _{x\in \dom (f)}|h_{f(x)}-h_x|<\infty .
  \]
  Then, by Minkowski's inequality, we have
  \begin{align*}
  \left (\sum _{x\in \dom (f) }|h_{f(x)}\xi _x|^2\right )^{1/2}\leq
  & \left (\sum _{x\in \dom (f) }|\xi _x|^2\cdot |h_{f(x)}-h_{x}|^2\right )^{1/2}\\
  &\quad + \left (\sum _{x\in \dom (f) }|h_x\xi _x|^2\right )^{1/2}\\
  &\leq s\|\xi \|+\| h\xi \|\\
  & < \infty .
  \end{align*} This shows that $v_f\xi $ is in the domain of $h$ and we are done.
  \end{proof}

The following result will later motivate our notion of a coarse operator.

  \begin{proposition}
  \label{ProphIsCoarseIFF}
  Let $X$ be u.l.f.\ metric space and $h\colon X\to \R $ be a map. The following are equivalent.
  \begin{enumerate}
    \item \label{Item1ProphIsCoarseIFF} The map $h\colon X\to \R $ is coarse.
    \item \label{Item2ProphIsCoarseIFF} The domain of $ h$, seen as an operator on $\ell _2(X)$, is invariant under partial translations and, for all $r>0$,
  \[
  \sup _f\|[h ,v_f]\|<\infty\ \text{ and } \ [h,v_f]\in \cstu(X),
  \]
  where the supremum is taken over the set of all partial $r$-translations $f$ of $X$.
  \item \label{Item2.5ProphIsCoarseIFF} The domain of $h$ is invariant under partial translations and
  \[
  \|[h ,v_f]
  \|<\infty
  \ \text{ and } \ [h,v_f]\in \cstu(X)\]
  for all partial translations $f$ of $X$.
  \end{enumerate}
Above we use the abuse of notation of identifying each $[h,v_f]$ in \eqref{Item2ProphIsCoarseIFF} and \eqref{Item2.5ProphIsCoarseIFF} with its unique bounded extension to $\ell_2(X)$.  \end{proposition}

  \begin{proof}
  \eqref{Item1ProphIsCoarseIFF}$\Rightarrow $\eqref{Item2ProphIsCoarseIFF}: Assume $h\colon X\to \R $ is coarse. By Proposition \ref{PropCoarseImpliesDomPreervation}, $\dom (h)$ is invariant under partial translations.  Fix $r>0$. As $h\colon X\to \R $ is coarse, there is $s>0$ such that $|h_x-h_y|\leq s$ for all $x,y\in X$ with $d(x,y)\leq r$.  Hence, if $f$ is a partial $r$-translation of $X$, then $|h_{f(x)}-h_x|\leq s$ for all $x\in \dom (f)$.
  Since \[[h,v_f]\xi=(h_{f(x)}-h_x)v_f\xi\] for all $\xi\in c_{00}(X)$,  it follows that  $\propg([h,v_f])\leq r$ and 
  \begin{align*}
  \| [ h, v_f]
  \|
  \leq \sup _{x\in \dom (f)}|h_{f(x)}-h_x|\leq s.
  \end{align*}

\eqref{Item2ProphIsCoarseIFF}$\Rightarrow $\eqref{Item2.5ProphIsCoarseIFF}: This is immediate.

\eqref{Item2.5ProphIsCoarseIFF}$\Rightarrow $\eqref{Item1ProphIsCoarseIFF}: Suppose $h\colon X\to \R $ is not coarse. So, there is $r>0$ and sequences $(x_n)_n$ and $(y_n)_n$ in $X$
  such that $d(x_n,y_n)\leq r$ for all $n\in \N $ and
  \begin{equation}
  \label{Eq1ProphIsCoarseIFF}
  \lim _n|h_{x_n}-h_{y_n}|=\infty .
  \end{equation}
  As $X$ is u.l.f., going to subsequences if necessary, we can assume that $(x_n)_n$ and $(y_n)_n$ are both sequences of distinct elements of $X$. Therefore, letting
  \[
  A=\{x_n\ :\ n\in \N \}\ \text { and }\ B=\{y_n\ :\ n\in \N \},
  \]
  the map $f\colon A\to B$ given by $f(x_n)=y_n$ for all $n\in \N $, defines a partial translation; a partial $r$-translation to be precise. However, we have
  \[
  \|[h,v_f]
  \|\geq \|[ h,v_f]\delta _{x_n}\|=|h_{f(x_n)}-h_{x_n}|
  \]
  for all $n\in \N $. By \eqref{Eq1ProphIsCoarseIFF}, the latter goes to infinity. But this contradicts the hypothesis.
  \end{proof}

Before finishing this subsection, we now pin-down an important property satisfied by diagonal flows on uniform Roe algebras.

  \begin{proposition}
  \label{PropIfFlowIsGivenByFucntionFinPropOpAreInDomOfInfinitesimalGenerator} Let $(X,d)$ be a u.l.f.\ metric space and let $h\colon X\to \R $ be a map inducing the flow $\sigma _h$ on $\cstu (X)$. For each partial translation $f$ of $X$, the map
  \[
  t\in \R \mapsto \sigma _{h,t}(v_f)\in \cstu (X)
  \]
  is differentiable at zero with derivative equal to $i[h,v_f]$. Therefore the map
  \begin{equation}
  \label{Eq24Jan24}
  t\in \R \mapsto \sigma _{h,t}(a)\in \cstu (X)
  \end{equation} is also differentiable at zero for all $a$ in $\cstu [X]$, with derivative equal to $i[h,a]$.
  \end{proposition}

  \begin{proof}
  As $\sigma _h$ is a flow, $h$ must be coarse (Proposition \ref{Prop2.1BragaExel}). Let us first show that the first statement indeed implies the second. Since
  $\cstu [X]$ is the linear span of $av_f$, for $a\in \ell _\infty (X)$, and $f$ a partial translation of $X$, it is enough to notice that the map \eqref{Eq24Jan24} is differentiable at zero for elements of the form $av_f$. Finally, as $\sigma _{h,t}(av_f)=a\sigma _{h,t}(v_f)$, we see that the second statement indeed follows from the first one.

  Fix now a partial translation $f$ of $X$. An elementary calculation shows that
  \begin{equation}
  \label{Conjugacao}
  \sigma _{h,t}(v_f)\delta _x = e^{it (h(f(x))-h(x))}v_f\delta _x,
  \end{equation}
  if $x$ is in $\dom (f)$, while $\sigma _{h,t}(v_f)\delta _x=0$ otherwise.  Denoting by $h\circ f-h$ the function on $X$ given by
  \[
  (h\circ f-h)(x) = \left \{
  \begin{matrix}
    h(f(x))-h(x), &\text {if } x\in \dom (f), \cr
    0, &\text {otherwise, }\hfill
  \end{matrix}\right .
  \]
  if follows from \eqref{Conjugacao} that
  $\sigma _{h,t}(v_f) = v_fe^{it(h\circ f-h)}$.  Observing that $h\circ f-h$ is a bounded function, thanks to $h$ being coarse, we deduce that
  $\sigma _{h,t}(v_f)$ is differentiable as a function of the variable $t$, and that its derivative at zero coincides with $iv_f(h\circ f-h)$.  Another elementary calculation shows that
  \[
  v_f(h\circ f-h)= [h, v_f],
  \]
  from where the conclusion follows.
  \end{proof}

  \begin{remark}
    Proposition \ref{PropIfFlowIsGivenByFucntionFinPropOpAreInDomOfInfinitesimalGenerator} can be stated in terms of the infinitesimal generator of the flow $\sigma $. Indeed, the \emph{infinitesimal generator} of a flow $\sigma $ on a $\mathrm C^*$-algebra $A$ is the operator $\delta $ on $A$ given by
  \[
  \delta (a)=\lim _{t\to 0}\frac {\sigma _t(a)-a}{t}
  \]
  for all $a\in A$ for which the limit above exists. With this terminology and considering a flow $\sigma $ on a uniform Roe algebra $\cstu (X)$, Proposition \ref{PropIfFlowIsGivenByFucntionFinPropOpAreInDomOfInfinitesimalGenerator} says that $\cstu [X]$ is contained in the domain of $\delta $, and that $\delta (a)=i[h,a]$ for all $a$ in $\cstu [X]$. We refer the reader to \cite [Chapter 3]{BratteliRobinsonVol1} for more on that.
  \end{remark}

\section{Coarse operators on $\ell _2(X)$}\label{SectionCoarseOp}\setcounter{subsection}{1}

We will now abstract some of the contents of Section \ref{SectionFlowsGivenByMaps} to develop the theory of \emph{noncommutative coarse functions} or, more precisely, of \emph{coarse operators} on $\ell _2(X)$ for a given u.l.f.\ metric space $X$. The following is our main definition and it adapts Proposition \ref{ProphIsCoarseIFF} to the noncommutative case of self-adjoint operators on $\ell _2(X)$.

  \begin{definition}
  [Definition \ref{DefinitionCoarseOp}]
  Let $X$ be a u.l.f.\ metric space and $h$ be a closed operator on $\ell _2(X)$. We call $h$ a \emph{coarse operator} if $\dom (h)$ is invariant under partial translations, 
  \[
  \| [ h ,v_f]
  \|<\infty\ \text{ and }\ [h,v_f]\in \cstu(X)
  \]
  for all partial translations $f$ on $X$; as usual, we are identifying $[h,v_f]$ with its unique extenssion to $\ell_2(X)$.
  \end{definition}

We should stress that the above definition does not require $h$ to be admissible; however this property will instead be \emph{proved} later in Theorem \ref{Thmc00IsAGreatCore!}.

\begin{example}\label{ExampleBoundedOpAreCoarse} Operators in $\cstu(X)$ are obviously coarse operators.  \end{example}

  \begin{example}
  By Proposition \ref{ProphIsCoarseIFF}, a map $h\colon X\to \R $ is coarse if and only if its interpretation as an operator on $\ell _2(X)$ is a coarse operator. For this reason, self-adjoint coarse operators on $\ell _2(X)$ should be thought of as noncommutative versions of coarse maps $X\to \R $.
  \end{example}

Just as in Proposition \ref{ProphIsCoarseIFF}, a more quantitative version of coarseness holds for operators. Precisely:

  \begin{proposition}
  \label{PropCoarseOpIFF}
  Let $X$ be a u.l.f.\ metric space and $h$ be a closed operator on $\ell _2(X)$, the domain of which is invariant under partial translations. Then $h$ is a coarse operator if and only if, for all $r>0$,
  \[
  \sup _f\|[h ,v_f]\|<\infty ,
  \]
  where the supremum is taken over the set of all partial $r$-translations $f$ of $X$, and the unique extension of $[h,v_f]$ to a bounded operator on $\ell_2(X)$ is in $\cstu(X)$ for all such $f$ as well.
  \end{proposition}

  \begin{proof}
  The ``if'' part being immediate, we move on to the ``only if'' direction.  Assuming that $h$ is a coarse operator we therefore suppose, by way of contradiction, that the supremum in the statement is infinite.

\begin{claim}\label{Claim017Oct24}{The supremum in the statement is infinite, even if taken over the set of all partial $r$-translations $f$ of $X$, with \emph{finite domain}.}  \end{claim}

\begin{proof} To see this, we fix $M>0$, and we choose a partial $r$-translation $g$ (which may or may not have finite domain), such that $\| [h,v_g]\| >M$.  Since $[h,v_g]$ is bounded on $\dom (h)$, by hypothesis, and since $c_{00}(X)$ is dense in $\dom (h)$, we may pick a unit vector $\xi $ in $c_{00}(X)$, such that
  \begin{equation}
  \label{Changingfntogn}
  M < \| [h,v_g]\xi \| = \| hv_g\xi - v_gh\xi \| .
  \end{equation} Fixing $\varepsilon >0$, let $C$ be a finite subset of $X$ such that
  \[
  \xi =p_C\xi ,\quad \text {and}\quad \| h\xi -p_Ch\xi \| <\varepsilon .
  \]
  Letting $f$ be the restriction of $g$ to $\dom (g)\cap C$, we then have that $v_{f} = v_gp_C$.  Moreover
  \[
  v_g\xi = v_{f}\xi ,\quad \text {and}\quad \| v_gh\xi - v_{f}h\xi \| <\varepsilon .
  \]
  So, upon substituting $g$ for $f$, the right-hand-side of \eqref{Changingfntogn} changes by at most $\varepsilon $, so the inequality is preserved as long as we choose $\varepsilon $ sufficiently small.  \end{proof}

  Slightly changing tack, we now make the following:

\begin{claim}\label{Claim17Oct24} If $C\subseteq X$ is any finite set, then the supremum in the statement is \emph{finite}, if taken over the set of all $r$-partial translations $f$ with $\dom (f)\subseteq C$.\end{claim} \begin{proof} Indeed, given such an $r$-partial translation $f$, we have that 
\begin{align*}
  \| [h,v_f]\| & = \| hv_f-v_fh \| \\
  & \leq \| hv_f \| + \| v_fh \| \\
  & \leq \| hp_{B_r(C)} \|+ \| p_{C}h \|, \end{align*}
  the last inequality due to the fact that $\dom (f)\subseteq C$, and hence that $\ran (f)\subseteq B_{r}(C)$, the ball of radius $r$ around $C$, because $r$-partial translations are not allowed to move points more than $r$ apart. As $h$ is necessarily bounded on the finite-dimensional space $\ell _2(B_r(C))$, $\| hp_{B_r(C)} \|$ is finite.  In order to estimate $\| p_Ch \|$, notice that \begin{align*}\| p_Ch\|\leq \|p_Chp_C\|+\|p_Chp_{X\setminus C}\|\leq \|hp_{C}\|+ \|p_Chp_{X\setminus C}\|.\end{align*} As before, $\|hp_C\|$ is finite and, since $p_Chp_{X\setminus C}=[h,p_{X\setminus C}]p_{X\setminus C}$, the fact that $h$ is coarse implies that $\|p_Chp_{X\setminus C}\|$ is also finite. This completes the proof of the claim.  \end{proof}

Our next claim is a strengthening of Claim \ref{Claim17Oct24}.

\begin{claim}\label{Claim317Oct24}{If $C\subseteq X$ is any finite set, then the supremum in the statement is \emph{infinite}, even if taken over the set of all partial $r$-translations $f$ of $X$, with finite domain, and such that both the domain of $f$ and the range of $f$ (henceforth denoted $\ran (f)$) are contained in $X\setminus C$.}\end{claim}

  \begin{proof}
  Given $M>0$, let $N$ be the supremum referred to in Claim \ref{Claim17Oct24} relative to the finite set $B_r(C)$, and, using Claim \ref{Claim017Oct24}, choose some $r$-partial translation $g$, with finite domain, and such that $\| [h,v_g]\| >M+N$. Letting $g_1$ and $g_2$ denote the restrictions of $g$ to $\dom (g)\cap B_r(C)$, and $\dom (g)\setminus B_r(C)$, respectively, we see that $v_g=v_{g_1}+v_{g_2}$, so
  \[
  \| [h,v_{g_2}]\| = \| [h,v_g-v_{g_1}]\| \geq \| [h,v_g]\| - \| [h,v_{g_1}]\| > M + N - N = M,
  \]
  Recalling that $\dom (g_2)\subseteq X\setminus B_r(C)$, we have that  $\ran (g_2)\subseteq X\setminus C$, for the same reason observed above that $r$-partial translations are not allowed to move points more than $r$ apart.  \end{proof}

\begin{claim}{There is a sequence $(f_n)_n$ of partial $r$-translations of $X$ such that the domains of the $f_n$ are finite and pairwise disjoint, the same holding for their ranges, and such that $\lim _{n\to \infty }\| [h,v_{f_n}]\| =\infty $.}\end{claim}

\begin{proof} As the proof of this claim is an easy consequence of Claim \ref{Claim317Oct24}, we leave it for the reader.\end{proof}

From now on we fix a sequence $(f_n)_n$ as above.

\begin{claim}\label{Claim517Oct24}{For every finite subset $C\subseteq X$, one has that
  \[
  \lim _{n\to \infty }\| p_{X\setminus C}[h, v_{f_n}]p_{X\setminus C} \| =\infty .
  \]}\end{claim}

\begin{proof}
In order to verify this claim,
observe that $\ell _2(C)$ is a finite dimensional space, hence
$p_Ch$ and $hp_C$ are bounded operators on the domain of $h$, say with
  \[
  \| p_Ch\| , \ \| hp_C\| \leq M_C.
  \]
  Therefore, for every $n$, we have that 
  \[
  \| [h, v_{f_n}]p_C\| =
  \| h v_{f_n}p_C\| + \| v_{f_n}hp_C\| \leq M_{B_r(C)} + M_C,
  \]
  and similarly
  $\| p_C[h, v_{f_n}]\| \leq M_{B_r(C)} + M_C$. Therefore
  \begin{align*}
  \| p_{X\setminus C}[h, v_{f_n}] p_{X\setminus C}\| & =
  \| (1-p_C)[h, v_{f_n}](1-p_C)\| \\
  & \geq
  \| [h, v_{f_n}]\| -
  \| [h, v_{f_n}]p_C\| \\ &\quad -
  \| p_C[h, v_{f_n}]\| -
  \| p_C[h, v_{f_n}]p_C\| \\
  &\geq
  \| [h, v_{f_n}]\| - 3(M_{B_r(C)} + M_C), \end{align*}
  from where the claim follows.
  \end{proof}

\begin{claim}\label{Claim617Oct24}
  Upon passing to a subsequence of $(f_n)_n$, there is a sequence $(X_n)_n$ of finite and pairwise disjoint subsets of $X$, such that both $\dom (f_n)$ and $\ran (f_n)$ are contained in $X_n$, and
  \[
  \lim _{n\to \infty }\| p_{X_n}[h, v_{f_n}]p_{X_n}\| =\infty .
  \]
  \end{claim}

\begin{proof} The proof of this claim is an easy application of Claim \ref{Claim517Oct24}, so we leave it for the reader.  \end{proof}

In order to conclude the proof, we let $(f_n)_n$ be as in Claim \ref{Claim517Oct24}, and we consider the $r$-partial translation defined on $\bigsqcup _{n\in {\mathbb N}}\dom (f_n)$ to be be a common extension of all of the $f_n$.

Observing that, for $n\neq m$, one has that $\dom (f_m)\subseteq X_m$, and that $X_m$ and $X_n$ are disjoint, it follows that $v_{f_m}p_{X_n}=0$, and similarly that $p_{X_n}v_{f_m}=0$, so
  \[
  p_{X_n}[h, v_{f}]p_{X_n} = p_{X_n}[h, v_{f_n}]p_{X_n},
  \]
  whence, for every $n$,
  \[
  \| [h, v_{f}]\| \geq
  \| p_{X_n}[h, v_{f}]p_{X_n}\| =
  \| p_{X_n}[h, v_{f_n}]p_{X_n}\| .
  \]
  Observing that the right hand side goes to $\infty $, as $n\to \infty $, by Claim \ref{Claim617Oct24}, we get a contradiction with the assumption that $[h, v_{f}]$ is a bounded operator.
  \end{proof}

We must now understand what  the vital properties satisfied by coarse operators are and, on the other hand, which properties are enough to guarantee coarseness of a given self-adjoint operator on $\ell _2(X)$. The next proposition deals with the former. For that, we introduce the \emph{conditional expectation} for certain unbounded operators on $\ell _2(X)$.

  \begin{definition}
  \label{DefiCondExp}
  Consider a set $X$ and any operator $h$ on $\ell _2(X)$ such that $c_{00}(X)\subseteq \dom (h)$. We define a map $E(h)\colon X\to \R $ by letting
  \[
  E(h)_x=\langle h\delta _x,\delta _x\rangle \ \text {for all}\ x\in X.
  \]
  \end{definition}

We identify the function $E(h)$
  defined above with an operator on $\ell _2(X)$ as usual. Under this identification, it is clear that $E$, when restricted to the bounded operators on $\ell _2(X)$, is simply the canonical conditional expectation $\cB (\ell _2(X))\to \ell _\infty (X)$.

  \begin{proposition}
  \label{PropCoarseOpFinPropMinusDiagIsBounded}
  Let $X$ be a u.l.f.\ metric space and $h$ be a coarse operator on $\ell _2(X)$. Then:
  \begin{enumerate}
         \item \label{PropCoarseOpFinPropMinusDiagIsBoundedItem1} The map $E(h)$
  given by Definition \ref{DefiCondExp} is coarse.
         \item \label{PropCoarseOpFinPropMinusDiagIsBoundedItem2} For each $A\subseteq X$, $p_{X\setminus A} hp_A$ is bounded and
  \[
  \sup \{\|p_{X\setminus A} hp_A\|\ :\ A\subseteq X\}<\infty .
  \]
  In particular, as an $X$-by-$X$ matrix, the coordinates of $h$ off the main diagonal are uniformly bounded.
  \end{enumerate}
  \end{proposition}

  \begin{proof}
  \eqref{PropCoarseOpFinPropMinusDiagIsBoundedItem1} Suppose $E(h)$ is not coarse. Then, there are $r>0$, and sequences $(x_n)_n$ and $(y_n)_n$ in $X$ such that
  \begin{equation}
  \label{Eq29Feb24.1}
  \lim _n|E(h)_{x_n}-E(h)_{y_n}|=\infty
  \end{equation}
  and $d(x_n,y_n)\leq r$ for all $n\in \N $. As $X$ is u.l.f., going to subsequences if necessary, we can assume that $(x_n)_n$ and $(y_n)_n$ are sequences of distinct elements. Let $A=\{x_n\ :\ n\in \N \}$, $B=\{y_n\ :\ n\in \N \}$, and $f\colon A\to B$ be the partial translation given by $f(x_n)=y_n$ for all $n\in \N $. Then,
  \begin{align*}
     \left \|[h,v_f]
  \right \|&\geq \left \|[h,v_f]\delta _{x_n}\right \|\\
  & = \left \|\sum _{z\in X} \langle h\delta _{y_n},\delta _z\rangle \delta _z - \langle v_fh\delta _{x_n},\delta _z\rangle \delta _z\right \|
  \\
  &\geq | \langle h\delta _{y_n},\delta _{y_n}\rangle - \langle h\delta _{x_n},\delta _{x_n}\rangle |\\
  &=|E(h)_{y_n}-E(h)_{x_n}|
     \end{align*}
  for all $n\in \N $. By \eqref{Eq29Feb24.1}, this shows that $ [h,v_f]$ is unbounded on $\dom (h)$, which contradicts the fact that $h$ is a coarse operator.

\eqref{PropCoarseOpFinPropMinusDiagIsBoundedItem2} As $h$ is a coarse operator, Proposition \ref{PropCoarseOpIFF} gives $s>0$, such that
  \[
  \|[h,p_A]
  \|\leq s
  \]
  for all $A\subseteq X$. Since $p_{X\setminus A}hp_A = [h, p_A]p_A$, the result follows.
  \end{proof}

We can now prove that coarse operators are automatically admissible, hence allowing
  us to understand how these operators act not only on $c_{00}(X)$ but on all of its domain.

\begin{theorem}\label{Thmc00IsAGreatCore!}
  Let $X$ be a u.l.f.\ metric space and $h$ be a coarse operator on $\ell _2(X)$.  If $\xi \in \dom (X)$, then
  \[h\xi =\sum _{x\in X}\langle \xi ,\delta _x\rangle h\delta _x.\]
  Consequently, $h$ is admissible. \end{theorem}

  \begin{proof}
  As $h$ is a coarse operator, Proposition \ref{PropCoarseOpFinPropMinusDiagIsBounded}\eqref{PropCoarseOpFinPropMinusDiagIsBoundedItem2} implies that
    \[M=\sup \{\|p_{X\setminus A}hp_A\|\ : \ A\subseteq X\}\]
    is finite.  Suppose towards a contradiction that the theorem fails and let $\xi \in \dom (h)$ be an offender. As $h$ is a closed operator, the series in the statement of the theorem not only fails to coincide with $h\xi $, but it does not even converge. Consequently that series also fails to satisfy Cauchy's criterion for the summability of a series, whence, letting $\xi _F=\sum _{x\in F}\langle \xi ,\delta _x\rangle \delta _x$ for each $F\subseteq X$, we can find $\eps >0$, and a sequence $(F(n))_n$ of finite disjoint subsets of $X$ such that
  \begin{equation}\label{EqPainInTheAss}
  \|h\xi _{F(n)}\|\geq \eps \ \text { for all }\ n\in \N .
  \end{equation}
  As $\lim _n\xi _{F(n)}=0$, going to a subsequence if necessary, we can assume that $\zeta =\sum _{n\in \N }\xi _{F(n)}$ has norm at most $\eps /(3M)$. In particular, \eqref{EqPainInTheAss} can be improved as follows:
  \begin{equation*}
  \|p_{F(n)}h\xi _{F(n)}\|\geq
  \|h\xi _{F(n)}\|-\|p_{X\setminus F(m)}hp_{F(n)}\xi _{F(n)}\|\geq
  \frac {2\eps }{3}
  \end{equation*} for all $n\in \N $. Since $h$ is a coarse operator, its domain is invariant under partial translations and, in particular, invariant under the projections $p_A$, $A\subseteq X$. So, $\zeta \in \dom (h)$, since $\zeta $ is the image of $\xi $ under such projection. However, \begin{align*}
    \|p_{F(n)}h\zeta \| & \geq \|p_{F(n)}h\xi _{F(n)}\|-\left \|p_{F(n)}hp_{X\setminus F(n)}\left (\sum _{m\in \N \setminus \{n\}}\xi _{F(m)}\right )\right \|\\
    &\geq \frac {\eps }{3} \end{align*} for all $n\in \N $. As $(F(n))_n$ are disjoint sets, this shows that $h\zeta $ has infinite norm, a contradiction.  \end{proof}

\begin{corollary}\label{CorCoincidenceOfDomains}
Let $X$ be a u.l.f.\ metric space, $h$ and $k$ be coarse operators on $\ell_2(X)$. If $h-k$ is bounded, then $\dom(h)=\dom(k)$.
\end{corollary}

\begin{proof}
 This is an immediate consequence of Proposition \ref{Prop05Nov24} and Theorem \ref{Thmc00IsAGreatCore!}.
\end{proof}

We finish this section with a corollary which will be fundamental in many of the forthcoming proofs.

\begin{corollary}
     \label{COR.PropCoarseOpFinPropMinusDiagIsBounded}
  Let $X$ be a u.l.f.\ metric space and $h$ be a coarse operator on $\ell _2(X)$. If there is an operator $h'$ on $\ell _2(X)$ with finite propagation, and such that $h-h' $ is bounded, then $h-E(h)$ is bounded. In particular, if $h$ has finite propagation, then $h-E(h) $ is bounded.  \end{corollary}

\begin{proof} Suppose there is an operator $h'$ on $\ell _2(X)$ with finite propagation so that $h-h'$ is bounded. In particular, $E(h-h')=E(h)-E(h')$ is also bounded.  Recall that, as $X$ is u.l.f., an operator $w$ on $\ell _2(X)$ with finite propagation is bounded if and only if $\sup _{x,y\in X}|\langle w\delta _x,\delta _y\rangle |<\infty $ (Proposition \ref{PropOperatorInducedFinPropULF}). By Proposition \ref{PropCoarseOpFinPropMinusDiagIsBounded}\eqref{PropCoarseOpFinPropMinusDiagIsBoundedItem2} above,
  \[
  \sup \left \{|\langle h\delta _x,\delta _y\rangle |\ :\ x,y\in X\text { with } x\neq y\right \}<\infty ,
  \]
  so, since $h-h'$ is bounded, we have
  \[
  \sup \left \{|\langle h'\delta _x,\delta _y\rangle |\ :\ x,y\in X\text { with } x\neq y\right \}<\infty .
  \]
  Hence,
  \[
  \sup _{x,y\in X}|\langle (h'-E(h'))\delta _x,\delta _y\rangle |<\infty .
  \]
  As $h'$ has finite propagation, so has $h'-E(h')$. So, $h'-E(h')$ is bounded. Finally, as the equality
  \[
  h-E(h)=h-h'+h'-E(h')+E(h')-E(h)
  \]
  holds on $c_{00}(X)$,  it follows that $h-E(h)$ is bounded on $c_{00}(X)$. Therefore, since both $h$ and $E(h)$ are coarse operators (for $E(h)$, this is a consequence of Propositions \ref{ProphIsCoarseIFF} and \ref{PropCoarseOpFinPropMinusDiagIsBounded}\eqref{PropCoarseOpFinPropMinusDiagIsBoundedItem2}), Theorem \ref{Thmc00IsAGreatCore!} implies that $h-E(h)$ is bounded on all of its domain. \end{proof}

\section{Flows given by self-adjoint operators on $\ell _2(X)$}\label{SectionFlowsGivenByOp}\setcounter{subsection}{1}

In this short section, we continue the study of Section \ref{SectionFlowsGivenByMaps} and extend it to the case of flows on uniform Roe algebras given by self-adjoint operators using the machinery of Section \ref{SectionCoarseOp}.  In view of Proposition \ref{PropIfFlowIsGivenByFucntionFinPropOpAreInDomOfInfinitesimalGenerator}, we start by introducing a geometric condition on our flows which will help us to better tame the self-adjoint operators inducing them.

  \begin{definition}
  [Definition \ref{DefinitionCoarseFlow}]
    Let $X$ be a u.l.f.\ metric space and $\sigma $ be a flow on $\cstu (X)$. We call $\sigma $ a \emph{coarse flow} if
  \[
  t\in \R \mapsto \sigma _{t}(v_f)\in \cstu (X)
  \]
  is differentiable for all partial translations $f$ of $X$.
  \end{definition}

  \begin{example}
    By Proposition \ref{PropIfFlowIsGivenByFucntionFinPropOpAreInDomOfInfinitesimalGenerator}, every flow on $\cstu (X)$ induced by a (necessarily coarse) map $X\to \R $ is coarse.
  \end{example}

Our next goal is to prove Proposition \ref{ProphCoarseOpImpliesFLow}, stated in the introduction, characterizing self-adjoint operators defining coarse flows on $\cstu (X)$. For that, we shall first need a lemma which is essentially  an adaptation of \cite [Proposition 3.2.55]{BratteliRobinsonVol1} to our settings. More precisely, we isolate below some consequences  which are contained in the proof of \cite [Proposition 3.2.55]{BratteliRobinsonVol1} but not explicitly stated therein. For the readers convenience, we add the proof of this subtleties below.

  \begin{lemma}\label{Lemma.Errata} Let $H $ be a Hilbert space, let $h$ be a
self-adjoint operator on $H$, and consider the pre-flow on $\cB (H)$ given by
letting
  \[
  \sigma _t(a)=e^{ith}ae^{-ith}
  \] for all $t\in \R $, and all $a\in \cB (H)$. For a given $a\in \cB (H)$, the
following are equivalent: \begin{enumerate}
  \item \label{Lemma.Errata.itemWOT} The map
  \[
  s\in \R \mapsto \sigma _s(a)\in \cB (H)
  \]
  is differentiable relative to the weak operator topology.
  \item \label{Lemma.Errata.item1} This map
  is differentiable relative to the strong operator topology.
  \item \label{Lemma.Errata.Item2} The domain of $h$ is invariant by $a$, and
$[h,a]$ is bounded on $\dom (h)$.
  \item \label{Lemma.Errata.Item3} In addition to the conditions in \eqref{Lemma.Errata.Item2},
  \[
  (\sigma _s(a)-a)\xi =\int _0^s\sigma _t(b) \xi dt,
  \]
  for all
 $s\in \R $ and all $\xi \in H$, where $b$ is the unique continuous extension of
$i[h,a]$ to $H$.
  \end{enumerate}
  In particular, if any one of the conditions above holds, then
  \[
  s\in \R \mapsto \sigma _s(a)\in \cB (H)
  \]
  is norm-continuous and its SOT-derivative at zero is $i[h,a]$.
  \end{lemma}

  \begin{proof}
  Before we begin, notice that, since
  \[
  \frac {\sigma _{s+t}(a)- \sigma _s(a)}t= \sigma _s\left (\frac {\sigma _t(a)- a}t\right ),
  \]
  the map in \eqref{Lemma.Errata.itemWOT} is differentiable at $s$ if and only if
it is differentiable at zero.

\medskip
    \eqref{Lemma.Errata.itemWOT}$\Rightarrow $\eqref{Lemma.Errata.Item2}\enspace
Observe that \eqref{Lemma.Errata.itemWOT} says that $a$ lies in the domain of the
infinitesimal generator of $\sigma $ (see \cite [Definition
3.1.5]{BratteliRobinsonVol1}).  Therefore, by \cite [Proposition
3.2.55]{BratteliRobinsonVol1}, there exists a core $D$ for $h$, such that $aD\subseteq D$,
and $[h,a]$ is bounded on $D$.

In order to prove that $\dom (h)$ is invariant by $a$, let $\xi \in \dom (h)$.  It then
follows that there exists a sequence $\{\xi _n\}_n\subseteq D$, converging to $\xi $, such that
$\{h(\xi _n)\}_n$ converges to $h(\xi )$.  We then have that
  \[
  h(a(\xi _n))= a(h(\xi _n)) + [h, a](\xi _n) \ \mathrel {\buildrel n\to \infty  \over \longrightarrow } \ a(h(\xi )) +
c(\xi ),
  \]
  where $c$ is the unique continuous extension of $[h,a]$ to $H$.
  This says that the sequence of ordered pairs
  \[
  \big \{\big (a(\xi _n), h(a(\xi _n))\big )\big \}_n
  \]
  converges to $\big (a(\xi ), a(h(\xi )) + c(\xi )\big )$. So, taking into account that $h$
is a closed operator, we deduce that $a(\xi )$ lies in $\dom (h)$, and
  \begin{align}
  \label{FormulaHA}
  h(a(\xi ))= a(h(\xi )) + c(\xi ).
  \end{align}
  This proves that $\dom (h)$ is invariant under $a$, so let us now prove that
$[h,a]$ is bounded on $\dom (h)$.  Given any $\xi $ in $\dom (h)$, observe that \eqref
{FormulaHA} gives
  $$
  [h,a](\xi ) = c(\xi ),
  $$
  so $[h,a]$ coincides with the bounded operator $c$ on $\dom (h)$, as required.

\medskip \eqref{Lemma.Errata.Item2}$\Rightarrow $\eqref{Lemma.Errata.Item3} Firstly,
as $i[h,a]$ is bounded on $\dom (h)$, let $b\in \cB (H)$ be the unique extension of
this operator to $H$.

For $\xi ,\zeta \in H$, a simple computation gives
  \begin{align}
  \label{variation}
  \begin{aligned}
  \big \langle (\sigma _{s+t}(a)-\sigma _{s}(a))\xi ,\zeta \big \rangle
  &=
  \big \langle e^{i(s+t)h}ae^{-i(s+t)h}\xi ,\zeta \big \rangle - \big \langle e^{ish}ae^{-ish}\xi ,\zeta \big \rangle
  \\ &=
  \big \langle ae^{-i(s+t)h}\xi ,(e^{-i(s+t)h}-e^{-ish})\zeta \big \rangle
  \\&+
  \big \langle a(e^{-i(s+t)h} -e^{-ish})\xi , e^{-ish}\zeta \big \rangle
  \\ &=
  \big \langle ae^{-i(s+t)h}\xi ,(e^{-ith}-1)e^{-ish}\zeta \big \rangle
  \\ &+
  \big \langle a (e^{-ith} -1)e^{-ish}\xi , e^{-ish}\zeta \big \rangle .
  \end{aligned}
  \end{align}

If moreover $\xi $ and $\zeta $ are in $\dom (h)$, then so are $e^{-ish}\xi $ and
$e^{-ish}\zeta $, whence
  \[
  \lim _{t\to 0}\frac {e^{-ith} -1}{t}e^{-ish}\xi =-ihe^{-ish}\xi ,
  \]
  and similarly after substituting $\zeta$ for $\xi$.
  Dividing both sides of \eqref{variation} by $t$, and taking
the limit as $t\to 0$, we get
  \begin{align*}
  \lim _{t\to 0} \left \langle
  \frac {\sigma _{s+t}(a)-\sigma _{s}(a)}{t}\xi ,\zeta \right \rangle &=
  \big \langle ae^{-ish}\xi , -ihe^{-ish}\zeta \big \rangle - \big \langle i ahe^{-ish} \xi ,
e^{-ish}\zeta \big \rangle \\ &=
  \big \langle ( iha -i ah)e^{-ish} \xi , e^{-ish}\zeta \big \rangle \\ &=
  \big \langle \sigma _s(b)\xi , \zeta \big \rangle .
  \end{align*}
  Integrating both sides then gives
  \begin{align*}
  \big \langle ( \sigma _{s}(a)- a )\xi ,\zeta \big \rangle =\int _0^s \big \langle \sigma _t(b)\xi ,
\zeta \big \rangle dt.
  \end{align*}
  As the assignment
  \[
  t\in \R \mapsto \sigma _t(b)\xi \in H
  \]
  is norm-continuous, the Bochner integral $\int _{0}^s\sigma _t(b)\xi dt$ is well
defined and hence
  \begin{align*}
  \langle ( \sigma _{s}(a)- a )\xi ,\zeta \rangle =\left \langle \int _0^s \sigma _t(b)\xi
dt, \zeta \right \rangle .
  \end{align*}
  Since the set of $\zeta $'s considered is dense in $H$, we get
  \begin{align*}
  (\sigma _{s}(a)- a)\xi = \int _0^s \sigma _t(b) \xi dt,
  \end{align*}
  for all $\xi \in \dom (h)$.  Again by density of $\dom (h)$, this holds for all $\xi $
in $H$.

\medskip \eqref{Lemma.Errata.Item3}$\Rightarrow $\eqref{Lemma.Errata.item1} \enspace
  Notice that, as
  \[
  t\in \R \mapsto \sigma _t(b)\xi \in H
  \]
 is continuous at $0$, dividing both sides of
  \[
  (\sigma _{s}(a)- a)\xi = \int _0^s \sigma _t(b) \xi dt
  \]
  by $s$, and letting it approach $0$, we get that
  \[
  \lim _{s\to 0}\frac {\sigma _{s}(a)- a}{s}\xi =b\xi .
  \]
 In other words, the map
  \[
  s\in \R \mapsto \sigma _{s}(a)\in \cstu (X)
  \]
 is SOT-differentiable at $0$ and its derivative is $b$. This proves
\eqref{Lemma.Errata.item1}.

\medskip \eqref{Lemma.Errata.item1}$\Rightarrow $\eqref{Lemma.Errata.itemWOT}
\enspace Obvious.

\medskip
 The last sentence in the statement follows since
  \[
  (\sigma _{s}(a)- a)\xi = \int _0^s \sigma _t(b) \xi dt,
  \]
  and, as $\sigma _t(b)$ is uniformly bounded in norm (by $\|b\|$), we must have
  \[
  \|(\sigma _{s}(a)- a)\xi \|\leq |s|\|b\|,
  \]
  for all $s\in \R $, and all $\xi \in H$ with norm at most $1$. Hence, the map
  \[
  s\in \R \mapsto \sigma _s(a)\in \cstu (X)
  \]
  is norm-continuous.
  \end{proof}

   \begin{proof}
  [Proof of Proposition \ref{ProphCoarseOpImpliesFLow}] 
  This is an immediate consequence of Lemma \ref{Lemma.Errata} applied to our settings. Precisely, suppose $\sigma $ is a coarse flow, i.e., $\sigma $ is a flow and
  \[
  t\in \R \mapsto \sigma _{t}(v_f)\in \cstu (X)
  \]
  is differentiable at zero for all partial translations $f$ on $X$. Since flows on $\cstu (X)$ are implemented by self-adjoint operators, we must have $\sigma =\sigma _h$ for some self-adjoint $h$ on $\ell _2(X)$ (Proposition \ref{ThmFlowsGivenByH}). The implication \eqref{Lemma.Errata.item1}$\Rightarrow $\eqref{Lemma.Errata.Item2} of Lemma  \ref{Lemma.Errata} then says that $\dom (h)$ is invariant under partial translations and that $[h,v_f]$ is bounded on $\dom (h)$. We are left to notice that $[h,v_f]$ is in $\cstu(X)$ for all partial translation $f$ on $X$, where here we identify $[h,v_f]$ with the unique bounded extension of $[h,v_f]$ to $\ell_2(X)$. This follows again by Lemma \ref{Lemma.Errata}. Indeed, the lemma and the fact that $t\in \R\mapsto \sigma_{h,t}(v_f)\in \cstu(X)$  is differentiable in norm at zero give that 
  \[i[h,v_f]=\SOTh\lim_{t\to 0}\frac{\sigma_{h,t}(v_f)-v_f}{t} =\lim_{t\to 0}\frac{\sigma_{h,t}(v_f)-v_f}{t}\in \cstu(X).\]

Conversely, suppose that $h$ is a coarse operator such that $\sigma_h$ defines a pre-flow on $\cstu(X)$ and fix a partial translation $f$ of $X$. Then $\dom (h)$ is $f$-invariant and $[h,v_f]$ is bounded on $\dom (h)$. We are left to notice that $t\in \R\mapsto \sigma_{h,t}(v_f)\in \cstu(X)$ is differentiable at zero. Notice that, since $h$ is a coarse operator, $[h,v_f]$ is in $\cstu(X)$, where here we identify $[h,v_f]$ with its unique bounded extension to $\ell_2(X)$. By  Lemma  \ref{Lemma.Errata}  the map 
\[s\in \R\mapsto \sigma_s([h,v_f])\in \cstu(X)\]
is  continuous. In particular, this map is Bochner integrable and  Lemma  \ref{Lemma.Errata}   gives that 
\[
  \sigma _{h,t}(v_f)-v_f =i\int _0^t\sigma _{h,t}([h,v_f]) dt.
  \]
  Therefore, 
\begin{align*}
\left\|\frac{\sigma_{h,t}(v_f)-v_f}{t}-i[h,v_f]\right\|   & =\left\|\frac{1}{t}\int_0^t\sigma_{h,s}([h,v_f])dt-[h,v_f]\right\|\\
&\leq \sup_{|s|\leq t}\left\|\sigma_{h,s}([h,v_f])-[h,v_f])\right\|\\
&\underset{t\to 0}{\longrightarrow} 0\end{align*}
and we we are done.
  \end{proof}

\begin{remark}
Since all elements of $\cstu(X)$ are coarse operators (Example \ref{ExampleBoundedOpAreCoarse}), Proposition \ref{ProphCoarseOpImpliesFLow} ensures that $\sigma_h$ is a coarse flow on $\cstu(X)$ for every $h \in \cstu(X)$. It is worth noting, however, that by choosing $h$ appropriately one may obtain flows that do not leave the algebraic subalgebra $\cstu[X]$ of $\cstu(X)$ invariant. Indeed, if $X$ is infinite and $x \in X$, one can choose $h$ such that $\sigma_{h,1}$ maps the projection onto $\mathrm{span}\{\delta_x\}$ to the projection onto $\mathrm{span}\{\xi\}$ for some $\xi$ with infinite support. Clearly, this latter projection does not belong to $\cstu[X]$. 
\end{remark}

We isolate a simple fact below for further reference.

  \begin{corollary}
  \label{CorValueDerFlow}
  Let $X$ be a u.l.f.\ metric space and $h$ be a self-adjoint operator on $\ell _2(X)$ such that $\sigma _h$ is a coarse flow on $\cstu (X)$. Then, for all partial translations $f$ of $X$, the derivative of
  \[
  t\in \R \mapsto \sigma _{h,t}(v_f)\in \cstu(X)
  \]
  at zero is $i[h,v_f]$.
  \end{corollary}

  \begin{proof}
    As in the proof of Proposition \ref{ProphCoarseOpImpliesFLow}, this is given by \cite [Proposition 3.2.55]{BratteliRobinsonVol1}.
  \end{proof}

  We finish this section showing that one can always find flows which are not coarse (as long as $X$ is infinite).

  \begin{proposition}
  \label{PropositionNonCoarseFLows}
  Let $X$ be a u.l.f.\ metric space with infinitely many points. Then there exists a flow on $\cstu (X)$ which is not coarse.
  \end{proposition}

  \begin{proof}
  Fix $x_0\in X$ and let $h(x)=d(x,x_0)$ for all $x\in X$. So, $h\colon X\to \R $ is a coarse map and $\sigma _h$ is a diagonal flow on $\cstu (X$) (Proposition \ref{Prop2.1BragaExel}). Clearly, $ \sigma _{uhu*}$ is also a flow on $\cstu (X)$ for all unitaries $u$ in $\cstu (X)$. Let us show that, if the unitary $u\in \cstu (X)$ is properly chosen, then $\sigma _{uhu^*}$ is not coarse. For that, notice that, since $h$ is unbounded, we can pick a vector $\xi \in \ell _2(X)$ with $\|\xi \|=1$ such that $\xi \not \in \dom (h)$. Moreover, without loss of generality, $\xi $ can be picked so that $\|\delta _{x_0}-\xi \|$ is arbitrarily small. Let $p=p_{x_0}$
  and let $q $ be the rank $1$ projection onto the span of $ \xi $. We can then assume that $ \|p-q\|$ is smaller than $1$, so there is a unitary $u$ in $\cstu (X)$ such that $uqu^*=p$. Since $u\xi $ equals $\lambda \delta _{x_0}$ for some $\lambda $ in the unitary circle, replacing $u$ by $\bar \lambda u$, we can assume that $ u\xi =\delta _{x_0}$. Hence, since $\xi $ is not in the domain of $h$, it follows that $\delta _{x_0}$ is not in the domain of $h'$.

We can now conclude the proof. If $\sigma _{h'}$ were to be coarse, the operator $h'$ would be a coarse operator (Proposition \ref{ProphCoarseOpImpliesFLow}). However, coarse operators have domains which are invariant under partial translations and, in particular, contain $c_{00}(X)$. Since $\delta _{x_0}$ is not in $\dom (h')$, it follows that $\sigma _{h'}$ cannot be a coarse flow.
  \end{proof}

  \begin{remark}
  Notice that, even though Proposition \ref{PropositionNonCoarseFLows} shows that non-coarse flows exist, the examples constructed therein are still very ``tame''. Indeed, they are given by unitary conjugation of a (diagonal) coarse operator $h$ on $\ell _2(X)$. In particular, the analysis of such flows, including for instance of their KMS states, reduces to the scenario of coarse (diagonal) operators; which we know how to deal with (see \cite {BragaExel2023KMS}).
  \end{remark}

\section{Automorphisms in the connected component of $\mathrm {id}_{\cstu (X)}$}\label{SectionAutConnectedIdentity}\setcounter{subsection}{1}

Section \ref{SectionFlowsGivenByOp} has dealt with Question \ref{QuestionsIandII}\ref{Item1Intro}, namely the question of when a pre-flow on a uniform Roe algebra is an actual flow and, for coarse flows, this is completely characterized in terms of the coarseness of the self-adjoint operator $h$ inducing $\sigma $ (see Proposition \ref{ProphCoarseOpImpliesFLow}). In this section and the following one, we will focus on  Question \ref{QuestionsIandII}\ref{Item2Intro}, that is, the question of when a self-adjoint operator $h$ on $\ell _2(X)$ defines a pre-flow on $\cstu (X)$. For that, we will now take a step back and study not flows per se, but strongly continuous paths of operators on uniform Roe algebras.

  \begin{definition} \label{DefineInTheComponent}
  Let $X$ be a u.l.f.\ metric space and $\sigma $ be an automorphism of $\cstu (X)$. We say that $\sigma $ \emph{is in the connected component of the identity on $\cstu (X)$} if there is a path $(\sigma _t)_{t\in [0,1]}$ of automorphisms of $\cstu (X)$ which is \emph{strongly continuous} in the sense that
  \[
  t\in \R \mapsto \sigma _t(a)\in \cstu (X)
  \]
  is continuous for all $a\in \cstu (X)$, and such that $\sigma _0=\mathrm {id}_{\cstu (X)}$ and $\sigma _1=\sigma $.
  \end{definition}

If $X$ is a u.l.f.\ metric space, any automorphism $\sigma $ of $\cstu (X)$ is given by a unitary $u$ in $\cB (\ell _2(X))$ by conjugation, i.e., $\sigma =\mathrm {Ad}(u)$ (see \cite [Lemma 3.1]{SpakulaWillett2013AdvMath}). In this section, we show that, if the automorphism $\sigma $ is in the connected component of the identity on $\cstu (X)$, then $u$   belongs to $\cstu (X)$, as long as $X$ has property A.  Precisely, the following is the main result of this section:

  \begin{theorem}
  \label{ThmFlowsuRaSpaciallyImplemented}
Let $ X $ be a u.l.f.\ metric space with property A,   and let $\sigma$ be an automorphism of $\cstu(X)$  lying in the connected component of the identity on $\cstu(X)$. If $w$ is any unitary operator on $\ell_2(X)$ such that $\sigma = \mathrm{Ad}(w)$, then $w \in \cstu(X)$.\footnote{  Note that, as mentioned above, the existence of such  a $w$ is   guaranteed by \cite[Lemma 3.1]{SpakulaWillett2013AdvMath}.}
  \end{theorem}

  Most of our uses of Theorem \ref{ThmFlowsuRaSpaciallyImplemented} will be through Theorem  \ref{ThmCoroFlowsuRaSpaciallyImplemented}, which we now proof as an immediate corollary of  Theorem \ref{ThmFlowsuRaSpaciallyImplemented}.

\begin{proof}[Proof of Theorem \ref{ThmCoroFlowsuRaSpaciallyImplemented}] Observing that each $\sigma _{h,t}$ clearly lies in the connected component of the identity (Definition \ref{DefineInTheComponent}), and also that $\sigma _{h,t}=\mathrm {Ad}(e^{ith})$, it follows from  Theorem~\ref{ThmFlowsuRaSpaciallyImplemented} that $e^{ith}\in \cstu (X)$.  \end{proof}


Our method to obtain Theorem \ref{ThmFlowsuRaSpaciallyImplemented} heavily depends on techniques developed in \cite {WhiteWillett2017}. Precisely, the following is the theorem therein whose proof will be essential for our approach to flows on uniform Roe algebras.

  \begin{theorem}
  \emph{(}\cite [Theorem E]{WhiteWillett2017}\emph{).}
    Let $X$ be a u.l.f.\ metric space with property A and $\sigma $ be an automorphism of $\cstu (X)$. Then there is a unitary $w\in \cstu (X)$ such that $w^*\sigma (\ell _\infty (X))w=\ell _\infty (X)$. \label{ThmWhiteWillett}
  \end{theorem}

We now explain some of the specifics of how Theorem \ref{ThmWhiteWillett} is obtained in order to isolate a lemma which will guide us in our proof of Theorem \ref{ThmFlowsuRaSpaciallyImplemented}. For that, let $\sigma \colon \cstu (X)\to \cstu (X)$ be an automorphism and let $f\colon X\to X$ be a map such that
  \begin{equation}
  \label{EqfCoarseEmbIntro}
  \inf _{x\in X}\left \|\sigma (e_{xx})\delta _{f(x)}\right \|>0.
  \end{equation}
  It is well known in the theory of rigidity of uniform Roe algebras that any map satisfying \eqref{EqfCoarseEmbIntro} must be a coarse equivalence (this was first obtained in \cite [Proof of Theorem 4.1]{SpakulaWillett2013AdvMath}, for a more citable result, see \cite [Theorem 4.12]{BragaFarah2018Trans}).  Suppose now that there is a bijection $g\colon X\to X$ which is close to $f$.\footnote {It is not yet known if such bijection always exists. As of now, this is known to be the case if $X$ has property A (\cite [Corollary 6.13]{WhiteWillett2017}), if $X$ is non-amenable (\cite [Theorem 5.1]{WhiteWillett2017}), or if $X$ is an expander graph (\cite [Theorems 1.2 and 1.4]{BaudierBragaFarahVignatiWillett2023BijExp}). } In particular, $g$ is also a coarse equivalence and we can define a unitary $v_g\colon \ell _2(X)\to \ell _2(X)$ by letting
  \begin{equation}
  \label{EqFormulaUnitaryIntro}
  v_g\delta _x=\delta _{g(x)}\ \text { for all }\ x\in X.
  \end{equation}
  As shown in \cite [Lemma 3.1]{SpakulaWillett2013AdvMath}, the automorphism $\sigma $ is implemented by a unitary on $\ell _2(X)$, so fix a unitary $u\in \cB (\ell _2(X))$ such that $\sigma =\mathrm {Ad}(u)$. A simple inspection of the proof of Theorem \ref{ThmWhiteWillett} tell us that, if $X$ has property A, then
  \[
  w=uv_g^*\in \cstu (X)
  \]
  (this is detailed in \cite [Page 983]{WhiteWillett2017}).

We summarize the discussion above in a lemma:

  \begin{lemma} \label{LemmaWhiteWillett2017}
  \emph{(}\cite [Proof of Theoorem E]{WhiteWillett2017}\emph{).}
    Let $X$ be a u.l.f.\ metric space with property A, $\sigma $ be an automorphism of $\cstu (X)$ implemented by a unitary $u\in \cB (\ell _2(X))$, $f\colon X\to X$ be a map satisfying \eqref{EqfCoarseEmbIntro}, and $g\colon X\to X$ be a bijection close to $f$. Then, letting $v_g^*$ be the unitary given by \eqref{EqFormulaUnitaryIntro}, the unitary $w=uv_g$ belongs to $\cstu (X)$.\footnote {Notice that $w$ in Lemma \ref{LemmaWhiteWillettIntro} is precisely $w$ in Theorem \ref{ThmWhiteWillett}.}\label{LemmaWhiteWillettIntro}
  \end{lemma}

We can now clearly state our approach towards Theorem \ref{ThmFlowsuRaSpaciallyImplemented}. Precisely, the following is our main technical result.

  \begin{theorem}
  \label{ThmCoarseEquivCloseId}
  Let $X$ be a u.l.f.\ metric space and $\sigma $ be an automorphism of $\cstu (X)$ in the connected component of the identity on $\cstu (X)$.  Then, for all $t\in \R $, there is a coarse equivalence $f_t\colon X\to X$ such that $f_t$ is close to the identity and
  \[
  \inf _{x\in X}\left \|\sigma _t(e_{xx})\delta _{f_t(x)}\right \|>0.
  \]
  \end{theorem}

We start the proof of Theorem \ref{ThmCoarseEquivCloseId} with a proposition.

  \begin{proposition}
  \label{PropCoarseEquivCloseId}
  Let $X$ be a u.l.f.\ metric space, $\sigma $ be an automorphism of $\cstu (X)$, and $f\colon X\to X$ be a map such that
  \[
  \delta =\inf _{x\in X}\left \|\sigma (e_{xx})\delta _{f(x)}\right \|>0.
  \]
  Suppose that for all $A\subseteq X$, there is a unitary $u\in \cstu (X)$ such that $\sigma (p_A)=up_Au^*$. Then $f$ is close to the identity.
  \end{proposition}

\begin{remark}
Proposition \ref{PropCoarseEquivCloseId} can also be derived from \cite[Theorem~10.2.1]{2024arXiv240313624M}. Given the length and general scope of \cite{2024arXiv240313624M}, we have chosen to retain our own proof, which may be more accessible to the general reader. We are grateful to the referee for bringing this to our attention.
\end{remark}

Before proving Proposition \ref{PropCoarseEquivCloseId}, we must recall the notion of equi-approxi\-ma\-ble families in $\cstu (X)$ and a criterion for equi-approximability obtained in \cite {BragaFarah2018Trans}; this will also play an important role in Section \ref{SecionUnboundedPart}.

  \begin{definition}
  Let $X$ be a u.l.f.\ metric space.
  \begin{enumerate}
  \item Given $\eps ,r>0$, we say that $a\in \cstu (X)$ is \emph{$\eps $-$r$-approximable} if there is $b\in \cstu (X)$, with $\propg (b)\leq r$, such that $\|a-b\|\leq \eps $.  \item A subset $S\subseteq \cstu (X)$ is \emph{equi-approximable} if for all $\eps >0$, there is $r>0$, such that every $a\in S$ is $\eps $-$r$-approximable.
  \end{enumerate}
  \end{definition}

  \begin{lemma}
  \emph{(}\cite [Lemma 4.9]{BragaFarah2018Trans}\emph{)}.  Let $X$ be a u.l.f.\ metric space and $(a_n)_n$ be a family of operators such that $a_M=\SOTh \sum _{n\in M}a_n\in \cstu (X)$ for all $M\subseteq \N $. Then the family $(a_M)_{M\subseteq \N }$ is equi-approximable.\label{LemmaBragaFarah49}
  \end{lemma}

  \begin{proof}
  [Proof of Proposition \ref{PropCoarseEquivCloseId}]
   Suppose towards a contradiction that $f\colon X\to X$ is not close to the identity. Then there is a sequence $(x_n)_n$ in $X$ such that
  \[
  \lim _{n\in \N }d(x_n,f(x_n))=\infty .
  \]
  Let $A=\bigcup _{n\in \N }\{x_n\}$. By hypothesis, there is a unitary $u\in \cstu (X)$ such that $\sigma (p_A)=up_Au^*$. As automorphisms of uniform Roe algebras are compact preserving and strongly continuous (see  \cite[Lemma 3.1]{SpakulaWillett2013AdvMath}), going to a subsequence if necessary, we can assume that
  \[
  \|e_{x_nx_n}u^*\sigma (e_{x_kx_k})u\|\leq 2^{-k-1}
  \]
  for all distinct $k,n\in \N $. Therefore,
  \begin{align}
  \label{Eq12PorpCloseIdentity}
  \left \|e_{x_nx_n}u^*\sigma (e_{x_nx_n})u\right \|&\geq \left \|e_{x_nx_n}u^*\sigma (p_A)u\right \|- \sum _{k\neq n}\left \|e_{x_nx_n}u^*\sigma (e_{x_kx_k})u\right \| \\
  &\geq \left \|e_{x_nx_n}p_A\right \|-\frac {1}{2}\notag \\
  &\geq \frac {1}{2}\notag
  \end{align}
  for all $n\in \N $.  Since $e_{x_nx_n}$, $u^*\sigma (e_{x_nx_n})$, and $e_{f(x_n)f(x_n)}$ have rank $1$, we have that
  \begin{align*}\left \|e_{x_nx_n}u^*\sigma (e_{x_nx_n})e_{f(x_n)f(x_n)}\right \|&=\left \|e_{x_nx_n}u^*\sigma (e_{x_nx_n})\right \|\left \|u^*\sigma (e_{x_nx_n})e_{f(x_n)f(x_n)}\right \|\\
  &=\left \|e_{x_nx_n}u^*\sigma (e_{x_nx_n})u\right \|\left \|\sigma (e_{x_nx_n})e_{f(x_n)f(x_n)}\right \|
  \end{align*} for all $n\in \N $ --- this is a simple computation and it can be found, for instance, in \cite [Lemma 6.5]{BragaFarahVignati2018AdvMath}.  Therefore, it follows from \eqref{Eq12PorpCloseIdentity} and our choice of $f$ that
  \begin{align}
  \label{Eq12PorpCloseIdentity2}
  \left \|e_{x_nx_n}u^*\sigma (e_{x_nx_n})e_{f(x_n)f(x_n)}\right \| \geq \frac {\delta }{2}
  \end{align}
  for all $n\in \N $.

As $\sigma $ is an automorphism, we know from \cite [Theorem 4.4]{BragaFarah2018Trans} that there is $r>0$ such that each $\sigma (e_{x_nx_n})$ is $\delta /3$-$r$-approximable. Since $u\in \cstu (X)$, replacing $r$ by a larger number if necessary, we can suppose that each $u^*\sigma (e_{x_nx_n})$ is $\delta /3$-$r$-approximable.  However, since $\lim _{n\in \N }d(x_n,f(x_n))=\infty $, this contradicts
  \eqref{Eq12PorpCloseIdentity2}. So, $f$ must be close to the identity of $X$.
  \end{proof}

  \begin{proof}
  [Proof of Theorem \ref{ThmCoarseEquivCloseId}] Fix $t\in \R $.  By \cite [Corollary 3.3]{BaudierBragaFarahKhukhroVignatiWillett2021uRaRig}, there is a map $f_t\colon X\to X$ such that
  \[
  \inf _{x\in X}\left \|\sigma _t(e_{xx})\delta _{f_t(x)}\right \|>0.
  \]
  As $\sigma $ is in the connected component of the identity, there is a continuous path of projections in $\cstu(X)$ starting in $p_A$ and ending in $\sigma(p_A)$. Therefore, it follows from \cite [Lemma 7.5.4]{MurphyBook} that $\sigma $ satisfies the hypothesis of Proposition \ref{PropCoarseEquivCloseId}. The result is now an immediate consequence of Proposition \ref{PropCoarseEquivCloseId}.
  \end{proof}

  \begin{proof}
  [Proof of Theorem \ref{ThmFlowsuRaSpaciallyImplemented}] By Theorem \ref{ThmCoarseEquivCloseId}, there is a coarse equivalence $f\colon X\to X$ satisfying \eqref{EqfCoarseEmbIntro} which is close to the identity on $X$. Since the identity $X\to X$ is a bijection, applying Lemma \ref{LemmaWhiteWillettIntro}
   to $g=\mathrm {id}_X$, we obtain that the unitary $u$ implementing $\sigma $ must belong to $\cstu (X)$ as desired.
  \end{proof}

\section{The unbounded part of $\cstu (X)$ }\label{SecionUnboundedPart}\setcounter{subsection}{1}

If a flow on $\cstu (X)$ is given by a self-adjoint operator $h$ on $\ell _2(X)$, then $e^{ith}\in \cstu (X)$ for all $t\in \R $, by Theorem \ref{ThmCoroFlowsuRaSpaciallyImplemented}. It is then important to understand which self-adjoint operators on $\ell _2(X)$ have this property. For that, we introduce the \emph{unbounded part of $\cstu (X)$}.

  \begin{definition}
  [Definition \ref{DefinitionUnboundedPartofURA}] Let $X$ be a u.l.f.\ metric space and $h$ be an admissible operator on $\ell _2(X)$. We say that $h$ is in the \emph{unbounded part of $\cstu (X)$} if, for all $\eps >0$, there is an operator $h'$ on $\ell _2(X)$ with finite propagation and such that $h-h'$ is a bounded operator with norm at most $\eps $. We denote the set of all such operators by $\mathrm {C}^*_{u,unb}(X)$.
  \end{definition}

  \begin{example}
  Clearly $\cstu (X)=\cB (\ell _2(X))\cap \mathrm {C}_{u,unb}^*(X)$. Moreover, if $h$ is an operator on $\ell _2(X)$ with finite propagation, it is immediate that $h\in \mathrm {C}_{u,unb}^*(X)$. In particular, if $h\colon X\to \R $ is any map, then its interpretation as a diagonal operator on $\ell _2(X)$ is in $\mathrm {C}_{u,unb}^*(X)$.
  \end{example}

  \begin{theorem}
  \label{ThmCoarseOpInUmbPart}
  Let $X$ be a u.l.f.\ metric space and $h $ be a self-adjoint operator in $\mathrm {C}^*_{u,unb}(X)$ such that
  \begin{equation}
  \label{EqSupFiniteOffDiag.}
  \sup \{|\langle h\delta _y,\delta _x\rangle |\ :\ x,y\in X, \ x\neq y\}<\infty .
  \end{equation}
  Then $e^{ih}\in \cstql (X)$.  In particular, if $h$ is a coarse operator, then $e^{ih}\in \cstql (X)$.
  \end{theorem}

The proof of Theorem \ref{ThmCoarseOpInUmbPart} is based on a test to guarantee an operator to be in $\cstql (X)$ based on Higson functions which was developed in \cite {SpakulaZhang2020JFA}. We recall this concept now.

  \begin{definition}
    Let $X$ be a u.l.f.\ metric space and $f\colon X\to \C $ be a function. We say that $f$ is \emph{Higson} if, for all $\eps >0$ and all $r>0$, there is a finite set $F\subseteq X$ such that, for all $x,y\in X\setminus F$,
  \[
  d(x,y)\leq r \ \text { implies }\ |f_x-f_y|\leq \eps .
  \]
  \end{definition}

  \begin{proof}
  [Proof of Theorem \ref{ThmCoarseOpInUmbPart}] By \cite [Theorem 3.3]{SpakulaZhang2020JFA}, an operator $a\in \cB (\ell _2(X))$ is in $\cstql (X)$ if and only if $[a,f]$ is compact for all Higson functions $f\colon X\to \C $. Fix such $f$ and let us show that $[e^{ih},f]$ is compact. This is equivalent to showing that $e^{ih}fe^{-ih}-f$ is compact, so, if $\pi \colon \cstu (X)\to \cstu (X)/\cK (\ell _2(X))$ is the canonical quotient map, we must show that $\pi (e^{ih}fe^{-ih})=\pi (f)$. For that, we prove that the map
  \begin{equation}
  \label{Eq29Feb24.6}
  t\in \R \mapsto \pi (e^{ith}fe^{-ith})\in \cstu (X)/\cK (\ell _2(X))
  \end{equation}
  is constant. Since this map equals $\pi (f)$ for $t=0$, the result will follow.

  \begin{claim}
  \label{Claim29Feb24}
  If $h'$ is an operator on $\ell _2(X)$ with finite propagation and such that $h-h'$ is bounded, then $[h',f]$ is bounded and compact.
  \end{claim}

  \begin{proof}
  As $h-h'$ is bounded on $c_{00}(X)$ and $h$ satisfies \eqref{EqSupFiniteOffDiag.}, this implies that
  \[
  \sup \{\langle h'\delta _y,\delta _x\rangle \ :\ x,y\in X \text { with }x\neq y\}<\infty .
  \]
  As $h'$ has finite propagation and
  \begin{align}
  \label{Eq24Fev24}
  \langle [h',f]\delta _y,\delta _x\rangle =(f_y-f_x)\langle h'\delta _y,\delta _x\rangle 
  \end{align}
  for all $x,y\in X$, the $X$-by-$X$ matrix representation of $[h',f]$ has finite propagation and its coordinates are uniformly bounded. Since $X$ is u.l.f., Proposition \ref{PropOperatorInducedFinPropULF} implies that the $X$-by-$X$ matrix $\big [\langle [h',f]\delta _y,\delta _x\rangle \big ]_{x,y\in X}$ canonically induces a bounded operator on $\ell _2(X)$. Being a finite propagation operator, $h'$ is also assumed to be admissible.  Based on the fact that $f$ is a diagonal operator, one can easily prove that $c_{00}(X)$ is a core for $[h',f]$. Therefore, since $[h',f]$ and the bounded operator induced by its matrix in the canonical basis of $\ell_2(X)$ clearly coincide on $c_{00}(X)$,  it follows that $[h',f]$ is itself bounded.

  In order to show that $[h',f]$ is compact, it is enough to notice that
  \[
  \lim _{x,y\to \infty }\langle [h',f]\delta _y,\delta _x\rangle =0
  \]
  (here we use again that $[h',f]$ has finite propagation).  But this follows immediately from \eqref{Eq24Fev24} since $f$ is Higson and $\propg (h')<\infty $.
  \end{proof}

  \begin{claim}
  \label{Claim29Feb24.2} $[h,f]$ is bounded and compact.
  \end{claim}

  \begin{proof}
  By hypothesis, for every $\varepsilon >0$, there is an operator $h'$ with finite propagation such that $h-h'$ is bounded and $\|h-h'\|<e$. Since the equality
  \begin{equation}\label{Eq17Oct24}
  [h,f]=[h-h',f]+[h',f]
  \end{equation}
  holds on $c_{00}(X)$, we have that
  \[
  \big \|[h,f]\big \|_0\leq \big \|[h-h',f]\big \|_0+\big \|[h',f]\big \|_0 \leq 2\varepsilon \|f\| + \big \|[h',f]\big \|_0 <\infty ,
  \]
  where by $\|\cdot \|_0$ we mean the norm of the restriction of the relevant operator to $c_{00}(X)$.  Again, since $f$ is diagonal,   $c_{00}(X)$ is a core for $[h,f]$, so the above inequality implies that $[h,f]$ is bounded on its whole domain, as required.  In addition,
  \[
  \big \|[h,f]-[h',f]\big \|_0 = \big \|[h-h',f] \big \|_0\leq 2\varepsilon \|f\|,
  \]
  which means that $[h,f]$ may be approximated arbitrarily close by compact operators, therefore itself being compact.
  \end{proof}

By Claim \ref{Claim29Feb24.2}, $[h,f]$ is bounded. Therefore, \cite [Proposition 3.2.55]{BratteliRobinsonVol1} implies that the map
  \begin{equation}
  \label{Eq29Feb24.9}
  t\in \R \mapsto e^{ith}fe^{-ith}\in \cB (\ell _2(X))
  \end{equation}
  is differentiable at zero and its derivative equals $i[h,f]$.  Moreover, since $\sigma _h$ is a group action, this implies that \eqref{Eq29Feb24.9} is differentiable at all $t\in \R $ with derivative $\sigma _{h,t}(i[h,f])$. As $[h,f]$ is compact (Claim \ref{Claim29Feb24.2}) and as the canonical quotient map $\pi\colon \cstu(X)\to \cstu(X)/\cK(\ell_2(X))$  is linear, we conclude that \eqref{Eq29Feb24.6} is differentiable at all $t\in \R $ with derivative $0$, hence a constant function.  Consequently
  \[
  \pi (e^{ih}fe^{-ih}) = \pi (f),
  \]
  whence $[e^{ih}, f]$ is compact, proving the required criterion for $e^{ih}$ to belong to $\cstql (X)$.

For the last statement, notice that, by Proposition \ref{PropCoarseOpFinPropMinusDiagIsBounded}, any coarse operator satisfies \eqref{EqSupFiniteOffDiag.}.
  \end{proof}


  \begin{theorem}
  \label{ThmhCoarseIsInUnboundedPartOfURA}
  Let $X$ be a u.l.f.\ metric space and let $h$ be a coarse operator  inducing a pre-flow on $\cstu(X)$, i.e.,  
  \[
  e^{ith}\cstu (X)e^{-ith}=\cstu (X)
  \]
  for all $t\in {\mathbb R}$. Then $h\in \mathrm {C}^*_{u,unb}(X)$. In particular, for $E$ as in Definition \ref{DefiCondExp}, $h-E(h)$ is bounded.
  \end{theorem}

  The proof of Theorem \ref{ThmhCoarseIsInUnboundedPartOfURA} is based on averaging operators over amenable groups and it closely follows ideas in \cite {LorentzWillett2020}. Before presenting its proof, we start setting up the territory to use the methods therein. Given a set $X$, consider the group $G=\{-1,1\}^X$, where the multiplication on $G$ is given by multiplication coordinate-wise. Although $G$ is often equipped with the product topology, we shall view it as a discrere group. Being a group, $G$ has a canonical action on $\ell _\infty (G)$ by left-multiplication given by letting
  \[
  (g\cdot a)_h=a_{g^{-1}h}
  \]
  for all $a=(a_g)_{g\in G}\in \ell _\infty (G)$ and all $g,h\in G$. As $G$ is abelian, it is also amenable and we can fix an invariant state $\varphi \colon \ell _\infty (G)\to \C $, meaning that
  \[
  \varphi (g\cdot a)=\varphi (a)
  \]
  for all $a\in \ell _\infty (G)$, and all $g\in G$, where $g\cdot a$ refers to the natural action of $G$ on $\ell ^\infty (G)$.
  Any given state on $\ell ^\infty (G)$ may be described as the integral against a certain \emph{finitely additive measure} on $G$, but our use of this theory will be limited to adopting their notation
  \[
  \int _Ga_gd\varphi (g),
  \]
  to mean nothing more that $\varphi (a)$, for $a=(a_g)_{g\in G}$ in $\ell _\infty (G)$.  By definition, the underlying set of $G$, namely $\{-1,1\}^X$, consists of all functions from $X$ to $\{-1,1\}$, and hence $G$ is, technically speaking, a subset of $\ell ^\infty (X)$.  Moreover, the inclusion map turns out to be a unitary representation of $G$ in $\ell ^\infty (X)$, and hence also in $\cstu (X)$.  We denote it by
  \begin{align}
  \label{EqPiRepresentationGroup}
  \pi \colon G&\hookrightarrow \ell ^\infty (X)\subseteq \cstu (X).
  \end{align}

The next lemma is well-known and its proof is a straightforward computation, so we omit the details.

  \begin{lemma}
  \label{LemmaCoordinatesAverageUnboundedoperatorByUnitaries}
  Let $X$ be a set and $E$ be as in Definition \ref{DefiCondExp}. Let $h$ be an operator on $\ell _2(X)$ such that $c_{00}(X)\subseteq \dom (h)$. Then
  \[
  \langle E(h)\delta _y,\delta _x\rangle =\int _G\langle \pi (\varepsilon )^* h\pi (\varepsilon )\delta _y,\delta _x\rangle
  d\varphi (\varepsilon )
  \]
  for all $x,y\in X$.\qed
  \end{lemma}

Given a set $X$, $\cL _1(\ell _2(X))$ denotes the trace-class operators on $\ell _2(X)$ and $\mathrm {Tr}\colon \cL _1(\ell _2(X))\to \C $ denotes the canonical trace on $\cL _1(\ell _2(X))$. We identify $\cB (\ell _2(X))$ with the dual of $\cL _1(\ell _2(X))$ in the canonical way, i.e., the duality is given by
  \begin{equation}
  \label{EqDualTRaceClass}
  a(b)=\mathrm {Tr}(ab)
  \end{equation}
  for all $a\in \cB (\ell _2(X))$, and all $b\in \cL _1(\ell _2(X))$. For each $a=(a_g)_{g\in G}\in \ell _\infty (G, \cB (\ell _2(X)))$,\footnote {Here $\ell _\infty (G, \cB (\ell _2(X)))$ denotes the $\mathrm C^*$-algebra of bounded functions $G\to \cB (\ell _2(X))$ endowed with the supremum norm.} we can define a functional
  \[
  \psi _a\colon b\in \cL _1(\ell _2(X))\mapsto \varphi ((a_g(b))_{g\in G})\in \C .
  \]
  Under the isometric identification above of $\cB (\ell _2(X))$ with the dual of $\cL _1(\ell _2(X))$, we consider each $\psi _a$ as an element in $\cB (\ell _2(X))$. We can then define the operator
  \begin{equation}
  \label{EqFormulaPsi}
  \Psi \colon a\in \ell _\infty (G, \cB (\ell _2(X)))\mapsto \psi _a\in \cB (\ell _2(X)).
  \end{equation}
  As $\varphi $ is contractive, it is clear that $\Psi $ is also contractive.
  Inasmuch as $\varphi \big ((a_g)_g\big )$ may be thought of as an \emph{average} of the bounded family of \emph{scalars} $a_g$, one may think of
  $\Psi \big ((a_g)_g\big )$ as an \emph{average} of the bounded family of \emph{operators} $a_g$. The following lemma will be important below.

  \begin{lemma}
  \emph{(}\cite [Lemma 3.3]{LorentzWillett2020}\emph{)}.  Let $X$ be a u.l.f.\ metric space, $G=\{-1,1\}^X$, and
  \[
  \Psi \colon \ell _\infty (G, \cB (\ell _2(X)))\to \cB (\ell _2(X))
  \]
  be as in \eqref{EqFormulaPsi}.  If $r>0$, and $(a_{\varepsilon })_{{\varepsilon }\in G}\in \ell _\infty (G, \cstu (X))$ is such that $\propg (a_{\varepsilon })\leq r$ for all $\varepsilon\in G$, then, $\Psi ((a_{\varepsilon })_{{\varepsilon }\in G})$ has propagation at most $ r$.\label{LorentzWillett2020Lemma3.3}
  \end{lemma}

  \begin{proof}
  [Proof of Theorem \ref{ThmhCoarseIsInUnboundedPartOfURA}]
  Let $\sigma _h$ be the pre-flow defined by $h$, which is actually a coarse flow thanks to Proposition \ref{ProphCoarseOpImpliesFLow}. By Proposition \ref{PropCoarseOpIFF}, there is $s>0$ such that
  \begin{equation*}\|[h,p_A]
  \|\leq s
  \end{equation*} for all $A\subseteq X$.  Hence, letting $\pi $ be the inclusion of $G=\{-1,1\}^X$ into $\ell_\infty(X)$ (see   \eqref{EqPiRepresentationGroup}), it follows that the family $([h,\pi (\varepsilon )])_{\varepsilon \in G}$ is uniformly bounded and in $\cstu(X)$, i.e.,  
  \begin{equation}
  \label{EqBounded17Feb24}
  ([h,\pi (\varepsilon )])_{\varepsilon \in G}\in \ell _\infty (G,\cstu (X)).
  \end{equation}
  Notice that this family is also equi-approximable. For that, observe that
  \[
  \SOTh \sum _{x\in A}[h,p_{\{x\}}]=[h,p_A]\in \cstu (X)
  \]
  for all $A\subseteq X$. Lemma \ref{LemmaBragaFarah49} then implies that $([h,\pi (\varepsilon )])_{\varepsilon \in G}$ is equi-approximable.

Fix an arbitrary $\theta >0$, and let us show that there is an operator $h'$ on $\ell _2(X)$, with finite propagation, and such that $h-h'$ is a bounded operator with norm at most $\theta $. For that, using the equi-approximability of the family $([h,\pi (\varepsilon )])_{\varepsilon \in G}$, pick $r>0$, such that $[h,\pi (\varepsilon )]$ is $\theta $-$r$-approximable for all $\varepsilon \in G$. Since each $\pi (\varepsilon )$ has propagation zero, $\pi (\varepsilon )^*[h,\pi (\varepsilon )]$ is also $\theta $-$r$-approximable for all $\varepsilon \in G$. For each $\varepsilon \in G$, pick $b_{\varepsilon }\in \cstu (X)$ with propagation at most $r$ and such that
  \begin{equation}
  \label{EqRelationbaw}
  c_{\varepsilon }=\pi (\varepsilon )^*[h,\pi (\varepsilon )]-b_{\varepsilon }
  \end{equation}
  has norm at most $\theta $. By \eqref{EqBounded17Feb24}, $([h,\pi (\varepsilon )])_{\varepsilon \in G}$ is bounded, so $(b_{\varepsilon })_{\varepsilon \in G}$ is bounded as well. Then, letting $\Psi $ be as in \eqref{EqFormulaPsi}, we can define
  \[
  w=\Psi ((\pi (\varepsilon )^*[h,\pi (\varepsilon )])_{\varepsilon \in G}),\ b=\Psi ((b_{\varepsilon })_{\varepsilon \in G}),\ \text { and }\ c=\Psi ((c_{\varepsilon })_{\varepsilon \in G}).
  \]
  Since each $b_{\varepsilon }$ has propagation at most $r$, Lemma \ref{LorentzWillett2020Lemma3.3} gives that the same holds for $b$. As each $c_{\varepsilon }$ has norm at most $\theta $, so has $c$. Also, since we have \eqref{EqRelationbaw} for all $\varepsilon \in \{-1,1\}^X$, it follows that
  \begin{equation}
  \label{EqRelationbaw2}
  c=w-b.
  \end{equation}

  \begin{claim}
  The operator $w+h$ has propagation zero.
  \end{claim}

  \begin{proof}
  First notice that, as $w$ is bounded and $h$ is a coarse operator, $w+h$ is also a coarse operator.
  Let now $x,y\in X$ be distinct and let us compute $\langle w\delta _y,\delta _x\rangle $. Letting $e_{y,x}$ be the partial isometry of rank $1$ such that $e_{y,x}\delta _x=\delta _y$, we have
  \begin{align*} \langle w\delta _y,\delta _x\rangle =\mathrm {Tr}(we_{y,x})=w(e_{y,x}),
  \end{align*} where the last equality is given by \eqref{EqDualTRaceClass}.  Therefore, by the formula of $\Psi $, $\varphi $, and using the duality relation in \eqref{EqDualTRaceClass} again, we have
  \begin{align*} \langle w\delta _y,\delta _x\rangle &=\Psi ((\pi (\varepsilon )^*[h,\pi (\varepsilon )])_{\varepsilon \in G})(e_{y,x}) \\
  &=\int _{\{-1,1\}^X} \mathrm {Tr}(\pi (\varepsilon )^*[h,\pi (\varepsilon )]e_{y,x}) d\varphi (\varepsilon ).
  \end{align*} Notice that, even though $h$ may not be bounded, the fact that $c_{00}(X)\subseteq \dom (h)$ implies that $he_{y,x}$ is well defined and bounded. Therefore, we can write
  \[
  \pi (\varepsilon )^*[h,\pi (\varepsilon )]e_{y,x}=\pi (\varepsilon )^*h\pi (\varepsilon )e_{y,x}-he_{y,x}
  \]
  and it follows that
  \begin{align*} \langle w\delta _y,\delta _x\rangle &=\int _{\{-1,1\}^X} \mathrm {Tr}(\pi (\varepsilon )^*h\pi (\varepsilon )e_{y,x})-\mathrm {Tr}(he_{y,x} )d\varphi (\varepsilon )\\
  &=\int _{\{-1,1\}^X} \langle \pi (\varepsilon )^*h\pi (\varepsilon )\delta _y,\delta _x\rangle d\varphi (\varepsilon )-\langle h\delta _y,\delta _x\rangle .
  \end{align*} Since $x$ and $y$ are distinct and $\langle E(h)\delta _y,\delta _x\rangle =0$, Lemma \ref{LemmaCoordinatesAverageUnboundedoperatorByUnitaries} implies that the first term in the right-hand side above is zero. This shows that $\langle w\delta _y,\delta _x\rangle =-\langle h\delta _y,\delta _x\rangle $ and we conclude that
  \[
  \langle (w+h)\delta _y,\delta _x\rangle =0.
  \]
  Since $x,y\in X$ were arbitrary distinct elements of $X$, this shows that $w+h$ has propagation zero as desired.
  \end{proof}

We can finally define $h'$ by letting
  \[
  h'=w+h-b.
  \] First notice that, as $b$ is bounded, and $w+h$ is a coarse operator, $h'$ is also a coarse operator. Moreover, since $w+h$ and $b$ actually have finite propagation, so does $h'$. By \eqref{EqRelationbaw2}, we have that $h-h'=-c$. So, $h-h'$ is a bounded operator with norm at most $\theta $. As $\theta >0$ was arbitrary, this shows that $h\in \mathrm {C}^*_{u,unb}(X)$.

  For the second statement, notice that the first part of the theorem implies that $h$ satisfies the hypothesis of Corollary \ref{COR.PropCoarseOpFinPropMinusDiagIsBounded}. The result is then a consequence of this proposition.
  \end{proof}

  \begin{proof}
  [Proof of Theorem \ref{ThmCharactCoarseFlowsInTermsOp}] \eqref{ThmCharactCoarseFlowsInTermsOpItem1}$\Rightarrow $\eqref{ThmCharactCoarseFlowsInTermsOpItem2}:  By \cite [Theorem 3.3]{SpakulaZhang2020JFA}, since $X$ has property A, we have $\cstu (X)=\cstql (X)$. Hence, Theorem \ref{ThmCoarseOpInUmbPart} gives that $\sigma _{h,t}(a)\in \cstu (X)$ for all $a\in \cstu (X)$ and, by Proposition \ref{ProphCoarseOpImpliesFLow}, $\sigma _h$ is indeed a coarse flow.

\eqref{ThmCharactCoarseFlowsInTermsOpItem2}$\Rightarrow $\eqref{ThmCharactCoarseFlowsInTermsOpItem1}: By Proposition \ref{ProphCoarseOpImpliesFLow}, $h$ is a coarse operator and, by Theorem \ref{ThmhCoarseIsInUnboundedPartOfURA}, $h\in \mathrm C^*_{u,unb}(X)$.

\eqref{ThmCharactCoarseFlowsInTermsOpItem2}$\Rightarrow $\eqref{ThmCharactCoarseFlowsInTermsOpItem3}: By Theorem \ref{ThmhCoarseIsInUnboundedPartOfURA}, $h$ is in $\mathrm C^*_{u,unb}(X)$ and $h-E(h)$ is bounded. As $h$ is in $\mathrm C^*_{u,unb}(X)$, for each $\eps >0$, there is an operator $h'$ with finite propagation such that $h-h'$ has norm at most $\eps $. As $h$ is a coarse operator (Proposition \ref{ProphCoarseOpImpliesFLow}), the supremum of the off-diagonal coordinates of $h$ is finite, hence the same is valid for $h'$. In particular, as $h'-E(h')$ has finite propagation, it follows that $h'-E(h')$ is bounded (Proposition \ref{PropOperatorInducedFinPropULF}). Since
  \[
  \|h-E(h)-(h'- E(h'))\|\leq \|h-h'\|+\|E(h)-E(h')\|\leq 2\eps
  \]
  and $\eps >0$ was arbitrary, this shows that $h-E(h)$ is in $\cstu (X)$. As $h$ is a coarse operator, it follows from Proposition \ref{PropCoarseOpFinPropMinusDiagIsBounded} that $E(h)\colon X\to \R $ is a coarse map.

\eqref{ThmCharactCoarseFlowsInTermsOpItem3}$\Rightarrow $\eqref{ThmCharactCoarseFlowsInTermsOpItem1}: As $h=h-E(h)+E(h)$ and $h-E(h)$ is assumed to be a bounded operator, the domain of $h$ is the one of $E(h)$. Hence, as we are also assuming here that $E(h)\colon X\to \R$ is a coarse map,  it follows that the domain of $h$ is invariant under partial translations (see Proposition \ref{PropCoarseImpliesDomPreervation}). Moreover, since $h-E(h)$ is bounded and $E(h)$ is coarse, it follows from Proposition \ref{ProphIsCoarseIFF} that
  \[
  [h,v_f]=[h-E(h),v_f]+[E(h),v_f]
  \]
  is bounded for all partial translations $f$ of $X$. So, $h$ is a coarse operator. Moreover, as $h-E(h)$ is not only bounded but belongs to $\cstu (X)$, it is immediate that $h$ is in $\mathrm C^*_{u,unb}(X)$.
  \end{proof}

\section{Cocycle perturbation and equivalence of flows on uniform Roe algebras}\label{SectionCocycle}\setcounter{subsection}{1}

There are natural equivalence relations on the set of flows on a given $\mathrm C^*$-algebra: \emph{inner perturbation}, \emph{cocycle perturbation}, and \emph{cocycle conjugacy} (see Definition \ref{DefiCocycleConjugate} above). In this section, we investigate these notions in the context of uniform Roe algebras. More precisely, given a u.l.f.\ metric space $X$, we aim to better understand the relationship between self-adjoint operators $h$ and $k$ on $\ell _2(X)$ which originate ``equivalent'' flows $\sigma _h$ and $\sigma _k$ and vice-versa. In case of diagonal flows, this is completely described in Theorem \ref{ThmCocycleEquivFlowCoarseMaps}. For more general coarse operators $h$ on $\ell _2(X)$ generating coarse flows on $\cstu (X)$, we obtain that $\sigma _h$ is always a cocycle perturbation of a diagonal flow as long as $X$ has property A (Theorem \ref{ThmFinPropagationCocycleEquivDiagFlowPropA}).  Finally, Theorem \ref{ThmClassificationFlowsGen} characterizes when coarse operators give rise to cocycle conjugate flows for spaces with property A; in other words, this is a noncommutative version of Theorem \ref{ThmCocycleEquivFlowCoarseMaps} under the presence of property A.

We begin with a general fact about self-adjoint operators and the unitary groups they generate.

\begin{lemma} \label{LemaOnOPG} Let $H$ be a separable Hilbert space, let $h$ and $k$ be   self-adjoint operators on $H$, and consider the map
  \[
  w_{h,k}: t\in {\mathbb R}\mapsto e^{ith}e^{-itk }\in \cB (H).
  \] Then, regarding the following statements, each one implies the next:
  \begin{enumerate}
  [label=\textnormal {(\alph *)}]
  \item \label{LemaOnOPGa}
  $w_{h, k}$ is differentiable at zero in the norm topology of $\cB (H)$.
  \item \label{LemaOnOPGb}
  $w_{h, k}$ is strongly differentiable at zero\footnote {That is, the map $t\mapsto e^{ith}e^{-itk}\xi $ is differentiable at zero for every $\xi $ in $H$.}.
  \item \label{LemaOnOPGc}
  $\dom (h)=\dom (k)$, and $h-k$ is bounded on that common domain.
  \item \label{LemaOnOPGd}
  $H$ contains a dense subspace $D$, such that $D\subseteq \dom (h)\cap \dom (k)$, $h-k$ is bounded on $D$, and $D$ is invariant under $e^{itk}$ for all $t$ in ${\mathbb R}$.
  \item \label{LemaOnOPGe}
  $w_{h, k}$ is Lipschitz continuous.
  \end{enumerate} \end{lemma}

\begin{proof} \ref{LemaOnOPGa}$\Rightarrow $\ref{LemaOnOPGb}.  Obvious

\medskip \noindent \ref{LemaOnOPGb}$\Rightarrow $\ref{LemaOnOPGc}.  By hypothesis, the formula
  \[
  \Delta (\xi ) = \lim _{t\to 0}\frac {e^{ith}e^{-itk}\xi -\xi }{t}
  \]
  provides a well defined linear operator on $H$. Moreover, as the operators \[\Delta_t\colon \xi\in H\mapsto \frac {e^{ith}e^{-itk}\xi -\xi }{t}\in H\]
  are bounded for all $t>0$, and $(\Delta_{1/n})_n$ converges pointwise to $\Delta$, the Banach--Steinhauss Theorem implies that $\Delta$ is a bounded operator.  In addition, given $\xi \in \dom (k)$, we have for all $t\in {\mathbb R}$ that
  \begin{equation}
  \label{eitheitk}
  e^{ith}\xi -\xi =
  e^{ith}\big (\xi -e^{-itk}\xi \big )+e^{ith}e^{-itk}\xi - \xi .
  \end{equation}
Since $(e^{ith})_t$ are uniformly bounded in norm and since $e^{ith}\xi\to \xi$, as $t\to 0$,  for all $\xi\in H$, we have that   
\begin{align*} 
e^{ith}\frac{(\xi -e^{-itk}\xi )}{t}=e^{ith}ik\xi-e^{ith}\left(\frac{e^{-itk}\xi -\xi}{t}+itk\right)\underset{t\to 0}{\to} ik\xi.
\end{align*}
So, dividing \eqref{eitheitk} by $t$ throughout, and taking the limit of the above as $t\to 0$, we get
  \[
  \frac {de^{ith}\xi }{dt}\Big |_{t=0} =
  ik\xi + \Delta \xi .
  \]
  Based on \cite [Proposition 5.3.13]{Pede:Analysis}, the conclusions are twofold: firstly,
  the existence of the derivative in the left-hand side implies that $\xi \in \dom (h)$, whence $\dom (k)\subseteq \dom (h)$, and, secondly, we deduce that
  \[
  ih\xi =
  ik\xi + \Delta \xi .
  \]
  Since the hypothesis is symmetric relative to $h$ and $k$, we see that $\dom (k)=\dom (h)$, and moreover $h-k=i^{-1}\Delta $, on that common domain, proving that $h-k$ is bounded.

\medskip \noindent \ref{LemaOnOPGc}$\Rightarrow $\ref{LemaOnOPGd}.  It suffices to take $D=\dom (k)$.

\medskip \noindent \ref{LemaOnOPGd}$\Rightarrow $\ref{LemaOnOPGe}.  Given $t_0, s\in {\mathbb R}$, observe that
  \begin{align*}
  e^{i(t_0+s) h}e^{-i(t_0+s) k} - e^{it_0 h}e^{-it_0 k}
  & =
  e^{it_0 h}\big (
  e^{is h} e^{-is k} - 1
  \big )e^{-it_0 k} \\
  & =
  e^{it_0 h}\Big (
  e^{is h}\big ( e^{-is k} - 1\big ) +
  e^{is h} - 1
  \Big )e^{-it_0 k}.
  \end{align*} Dividing by $s$ throughout,
  and taking the limit as $s\to 0$, we get,
  \[
  \frac {de^{ith}e^{-itk}\xi }{dt}\Big |_{t=t_0} = ie^{it_0 h}(h-k)e^{-it_0 k}\xi ,
  \] for any
  $
  \xi \in D.
  $
  Incidentally, it might be worth observing that $D$ is invariant under $e^{-it_0 k}$, by hypothesis, or else the above expression might not make sense.

  Letting $M$ be the norm of the restriction of $h-k$ to $D$, we deduce from the intermediate value theorem that
  \[
  \left \|e^{ith}e^{-itk }\xi -e^{ish}e^{-isk }\xi \right \|
  \leq M\, |t-s|\, \|\xi \|,
  \]
  for every $t$ and $s$ in ${\mathbb R}$. As $D$ is dense in $H$, we see that the above
in fact holds for every $\xi $ in $H$, from where the conclusion follows.
\end{proof}

  The following result will be the basis for the proofs of Theorems \ref{ThmFinPropagationCocycleEquivDiagFlowPropA}, \ref{ThmCocycleEquivFlowCoarseMaps}, and \ref{ThmClassificationFlowsGenMaisFraco}. Recall, a \cstar-algebra $A$ represented in some $\cB(H)$ is \emph{irreducible} if $\{0\}$ and $H$ are the only  subspaces of $H$ which are invariant with respect to all elements in $A$.

\begin{proposition}
  \label{BigPropCocycleEquiv}
  Let $H$ be a separable Hilbert space and let $A\subseteq \cB (H)$ be a unital $\mathrm C^*$-algebra.  Given self-adjoint operators $h$ and $k$ on $H$, suppose that
  \[
  e^{ith}A e^{-ith} = A, \quad \text {and} \quad e^{itk}A e^{-itk} = A
  \]
  for all $t\in {\mathbb R}$,
  and denote by $\sigma _h$ and $\sigma _k $ the associated pre-flows on $A$ defined by $h$ and $k$,  respectively.  Also consider the map
  \[
  w_{h,k}: t\in {\mathbb R}\mapsto e^{ith}e^{-itk }\in \cB (H).
  \]
  Then the following statements hold:
  \begin{enumerate} [label=\textnormal {(\alph *)}]
  \item \label{BigPropCocycleEquiv-a}
  If $A$ is irreducible and $\sigma _h$ is a cocycle perturbation of $\sigma _k$, then $w_{h,k}$ is norm continuous.
  \item \label{BigPropCocycleEquiv-b}
  If $A$ is irreducible and $\sigma _h$ is an inner perturbation of $\sigma _k$,
  then $\dom (h)=\dom (k)$ and $h-k$ is bounded, provided there is at least one nonzero vector\footnote {Please see Proposition \ref{EqualDomainWhenContainsCompacts} for a special case in which this somewhat annoying hypothesis is not needed.} belonging to both $\dom (h)$ and $\dom (k)$.
  \end{enumerate}
  Assuming that $w_{h,k}(t)$ lies in $A$ for all $t$ in $\R $, we also have the following:
  \begin{enumerate}
  [label=\textnormal {(\alph *)}]
  \setcounter {enumi}{2}
  \item \label{BigPropCocycleEquiv-c}
  If $w_{h,k}$ is norm continuous, then $\sigma _k $ is a cocycle perturbation of $\sigma _h$.
  \item \label{BigPropCocycleEquiv-d}
  If $w_{h,k}$ is differentiable at zero, then $\sigma _k $ is an inner perturbation of $\sigma _h$.
  \item \label{BigPropCocycleEquiv-e}
  If $H$ contains a dense subspace $D$, as in Lemma \ref{LemaOnOPG}\ref{LemaOnOPGd}, in particular such that $h-k$ is bounded on $D$, then $\sigma _k $ is a cocycle perturbation of $\sigma _h$.
  \item \label{BigPropCocycleEquiv-f}
  If $\dom (h)=\dom (k)$, and $h-k$ is bounded on this common domain, then $\sigma _k $ is a cocycle perturbation of $\sigma _h$.
  \item \label{BigPropCocycleEquiv-g}
  If $h$ and $k$ strongly commute\footnote {\label{CommutOps}
  We shall say that two  self-adjoint operators $h$ and $k$ strongly commute if there exists a projection valued measure $P$ on some measurable space $\Omega $, and two measurable, real valued functions $f$ and $g$ on $\Omega $, such that $h=\int _\Omega f(\omega )\,dP(\omega )$, and $k=\int _\Omega g(\omega )\,dP(\omega )$.  Very roughly, this means that $h$ and $k$ are ``simultaneously diagonalizable''.  In particular, note that $h$ and $k$ necessarily admit a common core, namely $\bigcup _{n\in {\mathbb N}} P(U_n)H$, where $U_n = \{\omega \in \Omega : |f(\omega )|+|g(\omega )|\leq n\}$.
  There are various alternative characterizations of strong commutativity, such as requiring that $e^{ith}$ commutes with $e^{isk}$, for every $t, s\in {\mathbb R}$, but nevertheless all such characterizations essentially express that there is an abelian von Neumann algebra to which both $h$ and $k$ are affiliated.}, and $H$ contains a dense subspace $D\subseteq \dom (h)\cap \dom (k)$ such that $h-k$ is bounded on $D$, then $\sigma _k $ is an inner perturbation of $\sigma _h$.
  \end{enumerate}
  \end{proposition}

  \begin{proof}
  \ref{BigPropCocycleEquiv-a}
  Let $(u_{t})_{t\in \R }$ be a cocycle for $\sigma _k$ such that
  \[
  \sigma _{h,t}=\mathrm {Ad}(u_t)\circ \sigma _{k,t}
  \]
  for all $t\in \R $. Unfolding definitions, this says that
  \[
  e^{ith}ae^{-ith}=u_te^{itk}ae^{-itk}u_t^*
  \]
  for all $a\in A$, and all $t\in \R $. Rewriting the above as
  \[
  ae^{-ith}u_te^{itk}=e^{-ith}u_te^{itk}a,
  \]
  we see that
  \begin{equation}
  \label{Eq09mar24}
  \lambda _t:= e^{-ith}u_te^{itk}
  \end{equation}
  lies in the commutant of $A$. Since $A$ is irreducible, we have that  $A'=\C \cdot \mathrm {id}_H$, so $\lambda _t$
  is a scalar multiple of the identity operator, clearly with absolute value one.
  As the right-hand-side of \eqref{Eq09mar24} is
  strongly continuous, $(\lambda _t)_{t\in \R }$ is continuous.  In addition, since
  \[
  w_{h, k} = e^{ith}e^{-itk}=\bar \lambda _tu_t
  \]
  for all $t\in \R $, and since $(u_t)_{t\in \R }$ is norm continuous, we conclude that $w_{h, k}$ is norm continuous as desired.

  \medskip \noindent \ref{BigPropCocycleEquiv-b}
  Let $u_t$ be a cocycle as in \ref{BigPropCocycleEquiv-a}, which we may now assume to be differentiable at zero.

We first observe that $u_0=1$, a fact that can easily be checked by plugging in $t=s=0$ in the cocycle condition in Definition \ref{DefiCocycleConjugate}(1).

As we are under a special case of \ref{BigPropCocycleEquiv-a}, we may re-run the argument given there up to \eqref{Eq09mar24} and then, picking any nonzero vector $\xi \in \dom (h)\cap \dom (k)$, we have that 
  \begin{align}
  \label{Eq10Mar24}
  \frac {d\lambda _t\xi }{dt}{\Big |_{t=0}} & =
  \lim _{t\to 0}\frac {e^{-ith}u_te^{itk}\xi -\xi }{t} \\
  &=\lim _{t\to 0}\Big ( e^{-ith}u_t\frac {e^{itk}\xi -\xi }{t} +e^{-ith}\frac {u_t\xi -\xi }{t} +\frac {e^{-ith}\xi -\xi }{t} \Big ) \notag \\
  & ={\frac {du_t\xi }{dt}}{\Big |_{t=0}} +i(k-h)\xi , \notag
  \end{align}
  from were we deduce that the map
  $t\in \R \mapsto \lambda _t\in \C $
  is differentiable at zero. Since
  \[
  w_{h, k}= e^{ith}e^{-itk}=\bar \lambda _tu_t,
  \]
  as already observed,
  this shows that $w_{h, k}$ is differentiable at zero.  So the conclusion follows from the implication \ref{LemaOnOPGa}$\Rightarrow $\ref{LemaOnOPGc} of Lemma \ref{LemaOnOPG}.


  \medskip \noindent \ref{BigPropCocycleEquiv-c}
  For each $t\in \R $, let $u_t=w_{h, k}(t)^*= e^{itk}e^{-i th}$. Then $u_t\in A$, for all $t\in \R $, and
  \begin{align*}
  u_t\sigma _{h,t}(u_s)&=e^{itk}e^{-i th}e^{ith}e^{isk}e^{-i sh}e^{-i th}\\
  &=e^{itk} e^{isk}e^{-i sh}e^{-i th}\\
  &=e^{i(t+s)k} e^{-i (t+s)h}\\
  &=u_{t+s}.
  \end{align*} Since $(u_t)_{t\in \R }$ is norm continuous, it is a cocycle for $\sigma _h$.  Hence, as $\sigma _{k,t}=u_t\sigma _{h,t}u^*_t$, for all $t\in \R $, $\sigma _k$ is a cocycle perturbation of $\sigma _h$.

  \medskip \noindent \ref{BigPropCocycleEquiv-d} Follows immediately, as above, under the extra assumption of differentiability of $w_{h, k}$ at zero.

  \medskip \noindent \ref{BigPropCocycleEquiv-e} Thanks to the implication \ref{LemaOnOPGd}$\Rightarrow $\ref{LemaOnOPGe} of Lemma \ref{LemaOnOPG}, we have that $w_{g,h}$ is continuous, so the conclusion follows from \ref{BigPropCocycleEquiv-c}, above.

  \medskip \noindent \ref{BigPropCocycleEquiv-f} This follows immediately from \ref{BigPropCocycleEquiv-e} upon taking $D=\dom (h)=\dom (k)$.

  \medskip \noindent \ref{BigPropCocycleEquiv-g}
  Choose any real number $M>\|(h-k){\restriction }_D\|$,
  and let $\Omega $, $P$, $f$, and $g$ be as in Footnote \ref{CommutOps}.  We then claim that $f-g$ is essentially bounded by $M$, in the sense that
  $P\big (\big \{\omega \in \Omega : |f(\omega )-g(\omega )|\geq M\big \}\big )=0$.  Assuming by contradiction that this is not so, and using the $\sigma $-additivity of $P$, we can find $N>0$, such that $P(U)\neq 0$, where
  \[
  U=\big \{\omega \in \Omega : |f(\omega )-g(\omega )|\geq M, \ |f(\omega )| \leq N, \text { and } |g(\omega )| \leq N \big \}.
  \]
  Choosing any unit vector in the range of $P(U)$, and taking a sequence of unit vectors $(\xi _n)_n$ in $D$ converging to $\xi $, we get
  \begin{align*}
  M & \leq
  \left \|\int _\Omega f(\omega )-g(\omega )\,dP(\omega ) \xi \right \| \\
  & = \|h\xi -k\xi \| \\
  & = \|hP_U\xi -kP_U\xi \| \\
  & = \|P_Uh\xi -P_Uk\xi \| \\
  & = \lim _{n\to \infty }\|P_Uh\xi _n-P_Uk\xi _n\| \\
  & \leq \big \|(h-k){\restriction }_D\big \|,
  \end{align*}
  where we have used that $P_Uh$ and $P_Uk$ are bounded operators (with norm at most $N$) in order to bring the above limit into play.
  This is a contradiction, hence proving that $f-g$ is essentially bounded.
  So $h-k$ is a bounded operator on its whole domain of definition and, since
  \[
  w_{h,k} = e^{ith}e^{-itk} = e^{it(h-k)},
  \]
  we deduce that $w_{h,k}$
  is differentiable everywhere, so the conclusion finally follows from \ref{BigPropCocycleEquiv-d}.  \end{proof}

Aiming at a proof of Theorem \ref{ThmClassificationFlowsGen}, let us now move to a slightly different situation, assuming, for the remainder of this section, that we are given
  a separable Hilbert space $H$ and a unital $\mathrm C^*$-algebra $A\subseteq \cB (H)$ containing the compacts.  Also let $h$ be a   self-adjoint operator on $H$ such that
  \[
  e^{ith}A e^{-ith} = A
  \]
  for all $t\in {\mathbb R}$,
  and denote by $\sigma _h$ the associated pre-flow on $A$.  Unlike in the previous situation, we will now assume that $\sigma _h$ is actually strongly continuous, that is, a flow on $A$.

We shall denote by $\delta _h$ the associated infinitesimal generator of $\sigma _h$, given by
  \[
  \delta _h(a) = \frac {d\sigma _{h, t}(a)}{dt}{\Big |_{t=0}},
  \]
  whose domain of course consists of those elements $a$ in $A$ for which the derivative exists.

Our next result shows that membership to $\dom (\delta _h)$ is closely related to membership to $\dom (h)$.

\begin{lemma} \label{OperatorAndDerivation} For every nonzero vector $\xi $ in $H$, a necessary and sufficient condition for $\xi $ to belong to $\dom (h)$ is that the orthogonal projection onto the one-dimensional space spanned by $\xi $, henceforth denoted $p_\xi $, lies in $\dom (\delta _h)$.  \end{lemma}

\begin{proof}
  By
  \cite [Proposition 3.2.55]{BratteliRobinsonVol1},
we have that  $p_\xi $ belongs to $\dom (\delta _h)$ if and only is there exists a core $D$ for $h$, such that $p_\xi (D)\subseteq D$, and $[h,p_\xi ]$ is a bounded operator on $D$.

Thus, assuming that $p_\xi \in \dom (\delta _h)$, and choosing a core $D$, as above, observe that $D$ cannot be contained in the kernel of $p_\xi $ because $D$ is dense in $H$.  So $p_\xi (D)$ is a nontrivial vector space contained in $D$, whence $\xi \in D\subseteq \dom (h)$.  Conversely, assuming that $\xi $ lies in $\dom (h)$, and taking the largest possible core $D=\dom (h)$, it is evident that $D$ is invariant under $p_\xi $, while, for every $\eta $ in $D$, we have  that
  \[
  [h,p_\xi ]\eta =
  hp_\xi \eta - p_\xi h\eta =
  h\big (\langle \eta ,\xi \rangle \xi \big ) - \langle h\eta ,\xi \rangle \xi =
  \langle \eta ,\xi \rangle h\xi - \langle \eta ,h\xi \rangle \xi ,
  \]
  which is evidently a bounded operator on $D$ as a function of the variable $\eta $.  So $p_\xi $ lies in $\dom (\delta _h)$ by the first paragraph in this proof.  \end{proof}

\begin{lemma} \label{EqualDomainWhenContainsCompacts}
  Let $H$ be a separable Hilbert space and let $A\subseteq \cB (H)$ be a unital $\mathrm C^*$-algebra containing the compacts. Also let $h$ and $k$ be self-adjoint operators on $H$ inducing flows $\sigma _h$ and $\sigma _k$ on $A$, as above.  If $\sigma _k$ is an inner-perturbation of $\sigma _h$, then $\dom (h)=\dom (k)$.  \end{lemma}

\begin{proof}
  By definition there exists a cocycle $(u_t)_{t\in {\mathbb R}}$ for $\sigma _h$, which is differentiable at zero, and such that
  \[
  \sigma _{k, t}(a)=u_t\sigma _{h,t}(a)u_t^*,
  \]
  for every $t\in {\mathbb R}$, and every $a\in A$.
  Letting $\delta _h$ and $\delta _k$ be the corresponding infinitesimal generators, the fact that $u_t$ is differentiable at zero clearly implies that $\dom (\delta _h)\subseteq \dom (\delta _k)$, so Lemma \ref{OperatorAndDerivation} implies that $\dom (h)\subseteq \dom (k)$.  Multiplying the equation displayed above by $u_t^*$ on the left, and by $u_t$ on the right, and applying a similar reasoning, then gives the reverse inclusion $\dom (h)\supseteq \dom (k)$.  \end{proof}


\subsection{Applications}

We now apply the machinery developed above to prove our main results, already stated in the introduction, classifying coarse flows on uniform Roe algebras of u.l.f metric spaces with property A. Our results are twofold: (1) we show that all coarse flows are cocycle perturbations of diagonal flows (Theorem \ref{ThmFinPropagationCocycleEquivDiagFlowPropA}) and (2) we characterize when two given flows are cocycle equivalent to each other in terms of the coarse operators generating them (Theorem \ref{ThmClassificationFlowsGen}). For diagonal flows, we also characterize both inner and cocycle perturbation (Theorem \ref{ThmCocycleEquivFlowCoarseMaps}).

  \begin{proof}
  [Proof of Theorem \ref{ThmClassificationFlowsGenMaisFraco}]
  Since $\sigma _h$ and $\sigma _k$ are flows, Theorem \ref{ThmFlowsuRaSpaciallyImplemented} guarantees that $e^{ith}$ and $e^{itk}$ lie in $\cstu (X)$ for all $t\in \R $, so
  \[
  w_{h,k}(t) = e^{ith}e^{-itk }\in \cstu (X)
  \]
  for all $t\in {\mathbb R}$. As both $h$ and $k$ are coarse operators and $h-k$ is bounded,  Corollary \ref{CorCoincidenceOfDomains} implies that $\dom(h)=\dom(k)$. Therefore,  Proposition \ref{BigPropCocycleEquiv}\ref{BigPropCocycleEquiv-f} implies that $\sigma_h$ and $\sigma_k$ are cocycle perturbations of each other.
  \end{proof}

We can now prove Theorem \ref{ThmCocycleEquivFlowCoarseMaps}, which summarizes our results for the commutative case, i.e., for diagonal flows on uniform Roe algebras.

  \begin{proof}
  [Proof of Theorem \ref{ThmCocycleEquivFlowCoarseMaps}] \eqref{ThmCocycleEquivFlowCoarseMapsItem1}$\Rightarrow $\eqref{ThmCocycleEquivFlowCoarseMapsItem2}: It is immediate that the diagonal operators $h$ and $k$, naturally associated to the given functions $h$ and $k$, strongly commute.  Moreover, given that the functions $h$ and $k$ are close, the difference operator $h-k$ is clearly bounded on $c_{00}(X)$.  In order to be able to apply Proposition 
  \ref{BigPropCocycleEquiv}\ref{BigPropCocycleEquiv-g},
  and hence to finish the proof, we need to prove that $w_{h,k}(t)$ lies in $\cstu (X)$, for every $t\in {\mathbb R}$, but this is obvious once we observe that the bounded operators $e^{ith}$ and $e^{itk}$ have propagation zero.

The implication \eqref{ThmCocycleEquivFlowCoarseMapsItem2} $\Rightarrow $\eqref{ThmCocycleEquivFlowCoarseMapsItem3} is obvious. For \eqref{ThmCocycleEquivFlowCoarseMapsItem3} $\Rightarrow $\eqref{ThmCocycleEquivFlowCoarseMapsItem1} notice that, since $\sigma _h$ and $\sigma _k$ are cocycle perturbations of each other, then Proposition \ref{BigPropCocycleEquiv}\ref{BigPropCocycleEquiv-a} implies that $w_{h, k}$ is norm continuous, and as we are in the commutative case, we have that 
  \[
  w_{h, k} = e^{ith}e^{-itk} = e^{it(h-k)}.
  \]
  In addition, by \cite [Theorem 13.36]{rudin1991functional}, a unitary group is norm continuous if and only if its generator, presently $h-k$, is a bounded operator, which clearly means that $h$ and $k$ are close as maps from $X$ to $\R $.
  \end{proof}

  \begin{proof}
  [Proof of Theorem \ref{ThmFinPropagationCocycleEquivDiagFlowPropA}]
  The conclusion will be reached upon applying Proposition \ref{BigPropCocycleEquiv}\ref{BigPropCocycleEquiv-e} to $h$ and $E(h)$, but first we must provide a dense subspace $D$ as required there.  For this we observe that $D:= c_{00}(X)$ is contained in the domain of $h$, because $h$ is coarse, and hence admissible by Theorem \ref{Thmc00IsAGreatCore!}.  $D$ is clearly also contained in the domain of the diagonal operator $E(h)$, and it is easy to see that the one-parameter unitary group generated by this operator leaves $D$ invariant.

Finally, observing that $h-E(h)$ is bounded on $D$ by Theorem \ref{ThmhCoarseIsInUnboundedPartOfURA}, we see that we are under the hypothesis of Proposition \ref{BigPropCocycleEquiv}\ref{BigPropCocycleEquiv-e}, from where the conclusion follows.
  \end{proof}


  \begin{proof}
  [Proof of Theorem \ref{ThmClassificationFlowsGen}] \eqref{ThmClassificationFlowsGenItem1}$\Rightarrow $\eqref{ThmClassificationFlowsGenItem2}: Suppose that $v\in \cstu (X)$ is a unitary element satisfying the conditions in Theorem 
  \ref{ThmClassificationFlowsGen}(\ref{ThmClassificationFlowsGenItem1}).  Letting $k'=vk v^*$, it is easy to see that
  \[
  \sigma _{k'}=\mathrm {Ad}(v)\circ \sigma _k\circ \mathrm {Ad}(v^*),
  \]
  so $\sigma _{k'}$ is clearly also a flow on $\cstu (X)$.  By Theorem \ref{ThmCoroFlowsuRaSpaciallyImplemented}, we then deduce that
  \[
  w_{h,k'}(t) := e^{ith}e^{-itk' }\in \cstu (X)
  \]
  for all $t\in {\mathbb R}$.
  So Proposition \ref{BigPropCocycleEquiv}\ref{BigPropCocycleEquiv-f} implies that $\sigma _{k'}$ is a \emph{cocycle perturbation} of $\sigma _h$, whence $\sigma _{k}$ is \emph{cocycle conjugate} to $\sigma _h$.

\eqref{ThmClassificationFlowsGenItem2}$\Rightarrow $\eqref{ThmClassificationFlowsGenItem1}: Suppose that $\sigma _h$ and $\sigma _k$ are cocycle conjugate to each other. By \cite [Lemma 1.1 and Corollary 1.3]{Kishimoto2000RepMathPhy}, this is equivalent to the existence of a unitary $v\in \cstu (X)$, and a family $(u_t)_{t\in \R }$ of unitaries in $\cstu (X)$ such that
  \begin{enumerate}
  \item $t\in \R \mapsto u_t\in \cstu (X)$ is differentiable everywhere, \item $(u_t)_{t\in \R }$ is a cocycle for $\mathrm {Ad}(v)\circ \sigma _{k}\circ \mathrm {Ad}(v^*)$, and \item \label{Item2Enumerate} $\sigma _{h,t}(a)=u_t v\sigma _{k,t} (v^*av)v^*u_t^*$ for all $t\in \R $ and all $a\in \cstu (X)$.
  \end{enumerate}
  Since the spectral theorem for self-adjoint operators gives
  \[
  e^{ivkv^*}=ve^{ik}v^*,
  \]
  item \eqref{Item2Enumerate} can be replaced by
  \begin{enumerate}
  \setcounter {enumi}{1} \item [3'.] $\sigma _{h,t}(a)=u_t \sigma _{vkv^*,t}(a)u_t^*$ for all $t\in \R $ and all $a\in \cstu (X)$.
  \end{enumerate}
  In other words, $\sigma _{h}$ is an inner perturbation of $\sigma _{vkv^*}$.  As $\cstu (X)$ contains the compacts, we may employ Lemma \ref{EqualDomainWhenContainsCompacts} to conclude that $\dom (h)=\dom (vkv^*)$, so the conclusion follows from Proposition \ref{BigPropCocycleEquiv}\ref{BigPropCocycleEquiv-b}.
  \end{proof}

\section{Pre-flows given by SOT-summable projections}\label{SectionPreflows}\setcounter{subsection}{1}

In this subsection, we identify a class of naturally occurring pre-flows on uniform Roe algebras which generalizes the diagonal pre-flows given by functions $X\to \R $ (Proposition \ref{PropFlowsGivenByProjections}). Our main result here provides an example of a pre-flow which is not a cocycle perturbation of a diagonal flow (Theorem \ref{ThmFlowExpGraph}), which goes on the opposite direction of Theorem \ref{ThmFinPropagationCocycleEquivDiagFlowPropA}.

  \begin{assumption}
  \label{Assumption1}
  Let $X$ be a u.l.f.\ metric space, $w\colon \N \to \R $ be a map, and $(p_n)_{n}$ be orthogonal projections on $\ell _2(X)$ such that
  \[
  p_M=\SOTh \sum _{n\in M}p_n\in \cstu (X)
  \]
  for all $M\subseteq \N $.  Let then, for each $n\in \N $,
  \[
  h_n=w(n)p_n \ \text { and }\ h=\SOTh \sum _{n} h_n.
  \]
  So, $h$ is a self-adjoint operator on $\ell _2(X)$.
  \end{assumption}

  \begin{example}
  Any pre-flow $\sigma $ given by a map $X\to \R $ is of the form $\sigma _h$ for some self-adjoint operator on $\ell _2(X)$ defined as above. Indeed, consider an arbitrary function $h\colon X\to \R $. For each $x\in X$, let $p_x=p_{\{x\}}$ and $h'_x=h(x)p_x$. It is clear that
  \[
  p_M=\SOTh \sum _{x\in M}p_x=p_{M}\in \cstu (X)
  \]
  for all $M\subseteq X$. Then, letting
  \[
  h'=\SOTh \sum _nh'_n,
  \]
  it is immediate that, as self-adjoint operators on $\ell _2(X)$, $h=h'$. So, $\sigma _h=\sigma _h'$ as desired.
  \end{example}

The previous example shows that a self-adjoint operator as in Assumption \ref{Assumption1} can be seen as a sort of generalized function $X\to \R $. In the beginning of Section \ref{SectionFlowsGivenByMaps}, we emphasize the (trivial) fact that if $h\colon X\to \R $ is an arbitrary map, then $\sigma _h$ is a pre-flow. The next proposition shows that the same holds for arbitrary self-adjoint operators $h$ as in Assumption \ref{Assumption1}.

  \begin{proposition}
  \label{PropFlowsGivenByProjections}
  In the setting of Assumption \ref{Assumption1}, $\sigma _h$ defines a pre-flow on $\cstu (X)$.
  \end{proposition}

  \begin{proof}
    This is a simple consequence of the ortogonality of the sequence of projections $(p_n)_n$. Indeed, this implies
  \[
  e^{ith}=\mathrm {id}_{\ell _2(X)}-p_\N +\SOTh \sum _{n}e^{ith_n}p_n,
  \]
  Moreover, since each $h_n$ is a bounded operator, their exponential can be easily computed and we have
  \begin{align}
  \label{EqehnFormula}
  e^{ith_n}p_n=\sum _{k=0}^\infty \frac {i^kt^kh_n^k}{k!}p_n=\sum _{k=0}^\infty \frac {i^kt^kw(n)^kp_n^k}{k!}p_n=p_n+(e^{itw(n)}-1)p_n
    \end{align}
  for all $n\in \N $ and all $t\in \R $. It then follows that
  \[
  e^{ith}=\mathrm {id}_{\ell _2(X)}+\SOTh \sum _n(e^{itw(n)}-1)p_n
  \]
  for all $t\in \R $ and the hypothesis on the strong convergence of sums of $(p_n)_n$ imply that $e^{ith}\in \cstu (X)$ for all $t\in \R $. We conclude that $\sigma _{h,t}(\cstu (X))\subseteq \cstu (X)$ for all $t\in \R $, i.e., $\sigma _h$ is a pre-flow on $\cstu (X)$.
  \end{proof}

The next result deals with expander graphs. We refer the reader to \cite[Chapter 5.6]{NowakYuBook} for precise definitions. In here, we simply say that if $(X_n)_n$ is a sequence of expander graphs and $X=\bigsqcup_nX_n$ is the coarse disjoint union of them, then 
\[p_M=\SOTh\sum_{n\in M}p_n\in \cstu(X)\]
for all $M\subseteq \N$, where for each $n\in\N$, $p_n$ is the rank-one projection of $\ell_2(X)$ onto the vector subspace of $\ell_2(X_n)\subseteq \ell_2(X)$ of vectors with constant coordinates; i.e., 
\[p_n=\begin{bmatrix}
    \frac{1}{|X_n|}& \cdots& \frac{1}{|X_n|}\\
    \vdots&\ddots&\vdots\\
    \frac{1}{|X_n|}&\cdots & \frac{1}{|X_n|}
\end{bmatrix}.\]
  \begin{theorem}
    \label{ThmFlowExpGraph}
    Let $X=\bigsqcup _n{X_n}$ be the coarse disjoint union of expander graphs and $w\colon \N \to \R $ be a map. For each $n\in \N $, let $p_n$ be the averaging projection on $\cB (\ell _2(X_n))$. Let $h_n=w(n)p_n$ for all $n\in \N $ and $h=\SOTh \sum _nh_n$. The following holds.
  \begin{enumerate}
        \item \label{ThmFlowExpGraphItem1} $\sigma _h$ is a pre-flow on $\cstu (X)$.
        \item \label{ThmFlowExpGraphItem1.5} If $w$ is unbounded, $\sigma _h$ is not a flow.
        \item \label{ThmFlowExpGraphItem2} If $w$ is unbounded, the pre-flow $\sigma _h$ is not a cocycle perturbation of any diagonal pre-flow on $\cstu (X)$.
          \end{enumerate}
  \end{theorem}

  \begin{proof}
  \eqref{ThmFlowExpGraphItem1} Since $X=\bigsqcup _nX_n$ is the coarse disjoint union of expander graphs, we know that $\sum _{n\in M}p_n$ is in $\cstu (X)$ for all $M\subseteq \N $ (see \cite {HigsonLafforgueSkandalis2002GAFA}). So, Proposition \ref{PropFlowsGivenByProjections} guarantees that $\sigma _h$ is a pre-flow on $\cstu (X)$.

\eqref{ThmFlowExpGraphItem1.5} Suppose now that $w$ is unbounded and let us notice $\sigma _h$ is not strongly continuous. For each $n\in \N $, pick $A_n\subseteq X_n$ such that $|A_n|=|X_n|/2$.\footnote {We are assuming here that each $X_n$ has even cardinality. If this is now the case, the argument will follow analogously by taking $|A_n|=\lceil |X_n|/2\rceil $ instead.}  Then, letting $A=\bigsqcup _n A_n$, we have
  \begin{align*}\left \|e^{ith}p_Ae^{-ith}-p_A\right \|&=\sup _n\left \|e^{ith_n}p_{A_n}e^{-ith_n}-p_{A_n}\right \|\\
  &=\sup _n\left \|e^{ith_n}p_{A_n}-p_{A_n}e^{ith_n}\right \|
    \end{align*} for all $t\in \R $. Since, for each $n\in \N $,
  \begin{align*}
    e^{ith_n}p_{X_n}=p_{X_n}+(e^{itw(n)}-1)p_n
    \end{align*}
    (cf.\ \eqref{EqehnFormula}), it follows that
  \begin{align*}
   \left \|e^{ith_n}p_{A_n}-vp_{A_n}e^{ith_n}\right \|=|e^{itw(n)}-1|\cdot \left \|p_np_{A_n}-p_{A_n}p_n\right \|
  \end{align*} for all $n\in \N $ and all $t\in \R $. In order to estimate the above, let $\mathbbm {1}_{A_n}$ denote the $|A_n|$-by-$|X_n{\setminus }A_n|$ matrix with $1$'s in all its entries; so, $\mathbbm {1}_{A_n}^*$ is the $|{X_n}{\setminus } A_n|$-by-$|A_n|$ matrix with $1$'s in all its entries. Hence, writing each $X_n$ in an order starting with elements in $A_n$ and later in $X_n{\setminus }A_n $, we can write
  \[
  (p_nv_{A_n}-v_{A_n}p_n)\restriction \ell _2(X_n)= \frac {1}{|X_n|}
  \begin{bmatrix}
    0&- \mathbbm {1}_{A_n}^*\\
  \mathbbm {1}_{A_n} & 0
  \end{bmatrix}.
  \]
  As the operator above has norm $1/2$, this shows that
  \[
  \left \|e^{ith_n}p_{A_n}-vp_{A_n}e^{ith_n}\right \|=\frac {1}{2}|e^{itw(n)}-1|
  \]
  for all $n\in \N $ and all $t\in \R $.

The conclusion of the previous paragraph is that
  \begin{align*}\left \|e^{ith}p_Ae^{-ith}-p_A\right \|
  =\frac {1}{2}\sup _n |e^{itw(n)}-1|
  \end{align*}
  for all $t\in \R $. As
  $w$ is unbounded, regardless of how small $\delta >0$ is, there is always $t\in \R $ with $|t|\leq \delta $ such that $|e^{itw(n)}-1|=2$. This shows that
  \[
  t\in \R \mapsto \sigma _h(p_A)\in \cstu (X)
  \]
  is not continuous; so, $\sigma _h$ is not a flow.

  \eqref{ThmFlowExpGraphItem2} By Proposition \ref{BigPropCocycleEquiv}\ref{BigPropCocycleEquiv-a}
  it is enough to show that if $k\colon X\to \R $ is a map, then $(e^{ith}e^{-itk})_{t\in \R }$ is not norm continuous. Fix such $k$ and for each $n\in \N $ let $k_n=k\restriction X_n$. Then, letting $(A_n)_n$ and $\mathbbm {1}_{A_n}$ be as above, we have
  \begin{align*}\left \|p_{A_n}(e^{ith_n}e^{-itg_n}-\mathrm {id}_{\ell _2(X_n)})p_{X{\setminus }A_n}\right \|&=|e^{itw(n)}-1|\cdot \left \|p_{A_n}p_np_{X_n{\setminus }A_n}\right \|\\
  &=\frac {|e^{itw(n)}-1|}{|X_n|}\cdot \left \|
  \begin{bmatrix}
          0& 0\\
  \mathbbm {1}_{A_n}&0
      \end{bmatrix}\right \|\\
  &=\frac {|e^{itw(n)}-1|}{2}
      \end{align*}
      for all $n\in \N $ and all $t\in \R $. Since,
  \[
  \left \|e^{ith}e^{-itk}-\mathrm {id}_{\ell _2(X)}\right \|=\sup _n\left \|e^{ith_n}e^{-itk_n}-\mathrm {id}_{\ell _2(X_n)}\right \|
  \]
  and $w$ is unbounded, this shows that $(e^{ith}e^{-itk})_{t\in \R }$ does not converge to $\mathrm {id}_{\ell _2(X)}$ in norm as $t\to 0$.
  \end{proof}

\begin{AIusage}
The first version of this manuscript,  published in \cite{BragaBussExel2026Flows}, has  a mistake in \cite[Proposition 1.3.7]{BragaBussExel2026Flows}. More precisely, \cite [Proposition 3.2.55]{BratteliRobinsonVol1} concerns the \emph{differentiability} of flows, and the mistake in our previous interpretation was in the precise meaning of the word ``differentiability'' in their context. In the published version of our article, we used differentiability with respect to the norm topology instead of the weak operator topology as it is actually the case in \cite [Proposition 3.2.55]{BratteliRobinsonVol1}. This misinterpretation led us to define what we called \emph{coarse operators} in a suboptimal way. Fortunately, this could be   easily fixed by redefining coarse operators as done in the current version of the manuscript.

While the published version the first version of this manuscript (see \cite{BragaBussExel2026Flows}) was written in the ``pre-AI era'' and the mistake in \cite[Proposition 1.3.7]{BragaBussExel2026Flows} was found by the authors with no AI assistance, GPT-6 Astra was used in the making of a corrigendum for \cite{BragaBussExel2026Flows}. More precisely, it was used to suggest possible corrections for \cite[Proposition 1.3.7]{BragaBussExel2026Flows}. This revised version was them made by modifying the original published version and merging it with the contents of the corrigendum.
\end{AIusage}

  \bibliographystyle{amsalpha}
  \bibliography{bibliography}

@ARTICLE{2024arXiv240313624M,
       author = {{Mart{\'\i}nez}, D. and {Vigolo}, F.},
        title = "{A rigidity framework for Roe-like algebras}",
      journal = {arXiv e-prints},
         year = 2024,
        month = mar,
          eid = {arXiv:2403.13624},
        pages = {arXiv:2403.13624},
          doi = {10.48550/arXiv.2403.13624},
archivePrefix = {arXiv},
       eprint = {2403.13624},
 primaryClass = {math.OA},
       adsurl = {https://ui.adsabs.harvard.edu/abs/2024arXiv240313624M}
}

@article {Kishimoto2000RepMathPhy,
    AUTHOR = {Kishimoto, A.},
     TITLE = {Locally representable one-parameter automorphism groups of
              {AF} algebras and {KMS} states},
   JOURNAL = {Rep. Math. Phys.},
  FJOURNAL = {Reports on Mathematical Physics},
    VOLUME = {45},
      YEAR = {2000},
    NUMBER = {3},
     PAGES = {333--356},
      ISSN = {0034-4877},
   MRCLASS = {82B10 (46L30 46L40 82B05 82B20)},
  MRNUMBER = {1772604},
MRREVIEWER = {Yoshikazu\ Katayama},
       DOI = {10.1016/S0034-4877(00)80003-8},
       URL = {https://doi.org/10.1016/S0034-4877(00)80003-8},
}

@book {MR4411373,
    AUTHOR = {Willett, R. and Yu, G.},
     TITLE = {Higher index theory},
    SERIES = {Cambridge Studies in Advanced Mathematics},
    VOLUME = {189},
 PUBLISHER = {Cambridge University Press, Cambridge},
      YEAR = {2020},
     PAGES = {xi+581},
      ISBN = {978-1-108-49106-8},
   MRCLASS = {46-02 (19Kxx 46L80 47L60 58J22)},
  MRNUMBER = {4411373},
MRREVIEWER = {Man-Ho\ Ho},
       DOI = {10.1017/9781108867351},
       URL = {https://doi.org/10.1017/9781108867351},
}

@ARTICLE{BaudierBragaFarahVignatiWillett2023BijExp,
       author = {{Baudier}, F. and {Braga}, B. M. and {Farah}, I. and {Vignati}, A. and {Willett}, R.},
        title = "{Coarse equivalence versus bijective coarse equivalence of expander graphs}",
      journal = {arXiv e-prints},
         year = 2023,
        month = jul,
          eid = {arXiv:2307.11529},
        pages = {arXiv:2307.11529},
          doi = {10.48550/arXiv.2307.11529},
archivePrefix = {arXiv},
       eprint = {2307.11529},
 primaryClass = {math.MG},
       adsurl = {https://ui.adsabs.harvard.edu/abs/2023arXiv230711529B}
}

@ARTICLE{BragaExel2023KMS,
       author = {{Braga}, B. M. and {Exel}, R.},
        title = "{KMS states on uniform Roe algebras}",
      journal = {arXiv e-prints},
         year = 2023,
        month = apr,
          eid = {arXiv:2304.05873},
        pages = {arXiv:2304.05873},
          doi = {10.48550/arXiv.2304.05873},
archivePrefix = {arXiv},
       eprint = {2304.05873},
 primaryClass = {math.OA},
       adsurl = {https://ui.adsabs.harvard.edu/abs/2023arXiv230405873D}
}

@book{BratteliRobinsonVol1,
    AUTHOR = {Bratteli, O. and Robinson, D.},
     TITLE = {Operator algebras and quantum statistical mechanics. 1},
    SERIES = {Texts and Monographs in Physics},
   EDITION = {Second},
      NOTE = {$C^\ast$- and $W^\ast$-algebras, symmetry groups,
              decomposition of states},
 PUBLISHER = {Springer-Verlag, New York},
      YEAR = {1987},
     PAGES = {xiv+505},
      ISBN = {0-387-17093-6},
   MRCLASS = {46Lxx (81-02 82-02)},
  MRNUMBER = {887100},
       DOI = {10.1007/978-3-662-02520-8},
       URL = {https://doi.org/10.1007/978-3-662-02520-8},
}

@article{BaudierBragaFarahKhukhroVignatiWillett2021uRaRig,
	author = {Baudier, F. and Braga, B. M. and Farah, I. and Khukhro, A. and {Vignati}, A. and Willett, R.},
	fjournal = {Inventiones Mathematicae},
	journal = {Invent. Math.},
	number = {3},
	pages = {1071--1100},
	title = {Uniform {R}oe algebras of uniformly locally finite metric spaces are rigid},
	volume = {230},
	year = {2022}}

@article{LorentzWillett2020,
	author = {Lorentz, M. and Willett, R.},
	journal = {Rocky Mountain Journal of Mathematics},
	number = {5},
	pages = {1747--1758},
	publisher = {Rocky Mountain Mathematics Consortium},
	title = {Bounded derivations on uniform {R}oe algebras},
	volume = {50},
	year = {2020}}

@article{Bourne2022JPhys,
	author = {Bourne, C.},
	fjournal = {Journal of Physics. A. Mathematical and Theoretical},
	journal = {J. Phys. A},
	number = {10},
	pages = {Paper No. 104004, 38},
	title = {Locally equivalent quasifree states and index theory},
	volume = {55},
	year = {2022}}

@article{LudewigThiang2021CommMathPhys,
	author = {Ludewig, M. and Thiang, G.},
	fjournal = {Communications in Mathematical Physics},
	journal = {Comm. Math. Phys.},
	number = {1},
	pages = {87--106},
	title = {Gaplessness of {L}andau {H}amiltonians on hyperbolic half-planes via coarse geometry},
	volume = {386},
	year = {2021}}

@article{Jones2021CommMathPhys,
	author = {Jones, C.},
	fjournal = {Communications in Mathematical Physics},
	journal = {Comm. Math. Phys.},
	number = {1},
	pages = {385--417},
	title = {Remarks on anomalous symmetries of {C}*-algebras},
	volume = {388},
	year = {2021}}

@book{BrownOzawa,
	author = {Brown, N. P. and Ozawa, N.},
	pages = {xvi+509},
	publisher = {American Mathematical Society, Providence, RI},
	series = {Graduate Studies in Mathematics},
	title = {{\cstar}-algebras and finite-dimensional approximations},
	volume = {88},
	year = {2008}}

@article{Roe1993,
	author = {Roe, J.},
	fjournal = {Memoirs of the American Mathematical Society},
	journal = {Mem. Amer. Math. Soc.},
	number = {497},
	pages = {x+90},
	title = {Coarse cohomology and index theory on complete {R}iemannian manifolds},
	volume = {104},
	year = {1993}}

@article{Roe1988,
	author = {Roe, J.},
	fjournal = {Journal of Differential Geometry},
	journal = {J. Differential Geom.},
	number = {1},
	pages = {87--113, 115--136},
	title = {An index theorem on open manifolds. {I}, {II}},
	volume = {27},
	year = {1988}}

@article{EwertMeyer2019,
	author = {Ewert, E. and Meyer, R.},
	fjournal = {Communications in Mathematical Physics},
	journal = {Comm. Math. Phys.},
	number = {3},
	pages = {1069--1098},
	title = {Coarse {G}eometry and {T}opological {P}hases},
	volume = {366},
	year = {2019}}

@article{Kubota2017,
	author = {Kubota, Y.},
	fjournal = {Communications in Mathematical Physics},
	journal = {Comm. Math. Phys.},
	number = {2},
	pages = {493--525},
	title = {Controlled topological phases and bulk-edge correspondence},
	volume = {349},
	year = {2017}}

@article{SpakulaZhang2020JFA,
	author = {{\v{S}}pakula, J. and Zhang, J.},
	fjournal = {Journal of Functional Analysis},
	journal = {J. Funct. Anal.},
	number = {1},
	pages = {108299, 25},
	title = {Quasi-locality and property {A}},
	volume = {278},
	year = {2020}}

@book{MurphyBook,
	author = {Murphy, G.},
	isbn = {0-12-511360-9},
	pages = {x+286},
	publisher = {Academic Press, Inc., Boston, MA},
	title = {{$C^*$}-algebras and operator theory},
	year = {1990}}

@book{NowakYuBook,
	author = {Nowak, P. and Yu, G.},
	doi = {10.4171/112},
	isbn = {978-3-03719-112-5},
	mrclass = {58-01 (19K56 20F69 46-01 58B34 58D10)},
	mrnumber = {2986138},
	mrreviewer = {Mikhail Ostrovskii},
	pages = {xiv+189},
	publisher = {European Mathematical Society (EMS), Z\"urich},
	series = {EMS Textbooks in Mathematics},
	title = {Large scale geometry},
	url = {https://doi.org/10.4171/112},
	year = {2012}}

@book{RoeBook,
	author = {Roe, J.},
	isbn = {0-8218-3332-4},
	pages = {viii+175},
	publisher = {American Mathematical Society, Providence, RI},
	series = {University Lecture Series},
	title = {Lectures on coarse geometry},
	volume = {31},
	year = {2003}}

@article{Yu2000,
	author = {Yu, G.},
	doi = {10.1007/s002229900032},
	fjournal = {Inventiones Mathematicae},
	journal = {Invent. Math.},
	number = {1},
	pages = {201--240},
	title = {The coarse {B}aum-{C}onnes conjecture for spaces which admit a uniform embedding into {H}ilbert space},
	volume = {139},
	year = {2000}}

@article{WhiteWillett2017,
	author = {{White}, S. and Willett, R.},
	journal = {Groups, Geometry, and Dynamics},
	pages = {949--989},
	title = {{Cartan subalgebras of uniform {R}oe algebras}},
	volume = {14},
	year = 2020}

@book{Pede:Analysis,
	address = {New York},
	author = {Pedersen, G. K.},
	publisher = {Springer-Verlag},
	series = {Graduate Texts in Mathematics},
	title = {Analysis now},
	volume = {118},
	year = {1989}}

@article{HigsonLafforgueSkandalis2002GAFA,
	author = {Higson, N. and Lafforgue, V. and Skandalis, G.},
	journal = {Geom. Funct. Anal.},
	number = {2},
	pages = {330--354},
	publisher = {Springer},
	title = {Counterexamples to the {B}aum--{C}onnes conjecture},
	volume = {12},
	year = {2002}}

@article{BragaFarah2018Trans,
	author = {Braga, B. M. and Farah, I.},
	fjournal = {Transactions of the American Mathematical Society},
	journal = {Trans. Amer. Math. Soc.},
	number = {2},
	pages = {1007--1040},
	title = {On the rigidity of uniform {R}oe algebras over uniformly locally finite coarse spaces},
	volume = {374},
	year = {2021}}

@article{SpakulaWillett2013AdvMath,
	author = {{\v{S}}pakula, J. and Willett, R.},
	doi = {10.1016/j.aim.2013.09.006},
	fjournal = {Advances in Mathematics},
	journal = {Adv. Math.},
	pages = {289--310},
	title = {On rigidity of {R}oe algebras},
	volume = {249},
	year = {2013}}

@article{BragaFarahVignati2018AdvMath,
	author = {Braga, B. M. and Farah, I. and Vignati, A.},
	doi = {10.1016/j.aim.2021.107886},
	fjournal = {Advances in Mathematics},
	journal = {Adv. Math.},
	pages = {Paper No. 107886, 35},
	title = {Uniform {R}oe coronas},
	volume = {389},
	year = {2021}}

@book{rudin1991functional,
	author = {Rudin, W.},
	publisher = {McGraw-Hill Science, Engineering \& Mathematics},
	title = {Functional analysis},
	year = {1991}}

@ARTICLE{Ozawa2023uRaSmallerQL,
       author = {{Ozawa}, N.},
        title = "{Embeddings of matrix algebras into uniform Roe algebras and quasi-local algebras}",
      journal = {arXiv e-prints},
         year = 2023,
        month = oct,
          eid = {arXiv:2310.03677},
        pages = {arXiv:2310.03677},
          doi = {10.48550/arXiv.2310.03677},
archivePrefix = {arXiv},
       eprint = {2310.03677},
 primaryClass = {math.OA},
       adsurl = {https://ui.adsabs.harvard.edu/abs/2023arXiv231003677O}
}

@ARTICLE{BragaBussExel2026Flows,
       author = {{Braga}, B. M. and {Buss}, A. and {Exel}, R.},
        title = "{Flows on uniform Roe algebras}",
      journal = {Journal of the Institute of Mathematics of Jussieu},
         year = 2026,
        pages = {1393–1439},
}

  \end{document}